\documentclass[12pt]{amsart}

\usepackage{amsmath}
\usepackage{latexsym}
\usepackage{amssymb}
\usepackage{amsfonts}

\usepackage{graphics}

\renewcommand{\proof}{\par\noindent{\it Proof.\ \ }}
\def\qed{\ifmmode\square\else\nolinebreak\hfill
$\square$\fi\par\vskip5pt}

\def\div{\,\big|\,}

 \def\ZZ{\mathbb Z}

\def\mod{{\sf mod~}}

\def\Aut{{\sf Aut}}  \def\Out{{\sf Out}}
\def\Cos{{\sf Cos}}

\def\K{{\sf K}}
\bf

\def\D{{\rm D}} \def\Q{{\rm Q}} 
 \def\S{{\rm S}} \def\G{{\rm G}}
\def\J{{\rm J}} \def\M{{\rm M}} \def\Fit{{\rm Fit}}
\def\soc{{\rm soc}} 
\def\N{{\bf N}} \def\O{{\bf O}}

\def\Ga{{\it\Gamma}}

\def\PG{{\rm PG}}

\def\POmega{{\rm P\Omega}}
\def\Sp{{\rm Sp}}\def\PSp{{\rm PSp}} 
\def\GammaL{{\rm \Gamma L}}

\def\A{{\rm A}}\def\Sym{{\rm Sym}}
\def\B{{\rm B}}
\def\PSL{{\rm PSL}}\def\PGL{{\rm PGL}}
\def\GL{{\rm GL}} \def\SL{{\rm SL}}
\def\L{{\rm L}}

 \def\PSU{{\rm PSU}} \def\SU{{\rm SU}}

\def\Sz{{\rm Sz}} \def\McL{{\rm McL}} \def\He{{\rm He}}
\def\Ru{{\rm Ru}} \def\Th{{\rm Th}}
\def\E{{\rm E}} \def\F{{\rm F}} \def\D{{\rm D}} \def\G{{\rm G}}
 
\def\Co{{\rm Co}} 
\def\HS{{\rm HS}}
\def\Ree{{\rm Ree}}
\def\ON{{\rm O'N}}
\def\Fi{{\rm Fi}}
\def\Ly{{\rm Ly}}

\newtheorem{theorem}{Theorem}[section]%
\newtheorem{lemma}[theorem]{Lemma}%
\newtheorem{corollary}[theorem]{Corollary}%

\newtheorem{hypothesis}[theorem]{Hypothesis}
\begin{document}

\title[edge-primitivity]{On edge-primitive graphs with soluble edge-stabilizers}
\thanks{Supported by the National Natural Science Foundation of China (11971248, 11731002) and the Fundamental Research Funds for the Central Universities.}
\author[Han]{Hua Han}
\address{H. Han\\ School of Science, Tianjin University of Technology\\
Tianjin 300384\\
P. R. China}
\email{hh1204@mail.nankai.edu.cn}

\author[Liao]{Hong Ci Liao}
\address{H. C. Liao\\ Center for Combinatorics\\
LPMC, Nankai University\\
Tianjin 300071\\
P. R. China} \email{827436562@qq.com}

\author[Lu]{Zai Ping Lu}
\address{
Zaiping Lu\\
Center for Combinatorics, LPMC,
Nankai University,
Tianjin 300071, China
}
\email{lu@nankai.edu.cn}

\begin{abstract}
A graph is edge-primitive if its automorphism group acts primitively on the edge set, and $2$-arc-transitive if its automorphism group acts transitively on the set of $2$-arcs. In this paper, we present a classification for those edge-primitive graphs which are $2$-arc-transitive and have soluble edge-stabilizers.

\vskip 2pt

\end{abstract}
\maketitle
\date\today

\parskip 5pt

\section{Introduction}

In this paper, all  graphs are assumed to be finite and simple,
and all groups are assumed to be finite.

A graph is a pair $\Ga=(V,E)$
of a nonempty set $V$ and a set $E$ of $2$-subsets of $V$.
The elements in $V$ and $E$ are called the vertices and edges of $\Ga$, respectively.
For $v\in V$, the set $\Ga(v)=\{u\in V\mid \{u,v\}\in E\}$ is  called
the neighborhood of $v$ in $\Ga$, while $|\Ga(v)|$ is called the valency
of $v$. We say that the graph $\Ga$ has valency $d$ or $\Ga$ is $d$-regular
if its vertices  have equal valency $d$.
 For an integer $s\ge 1$, an $s$-arc in $\Ga$ is
an $(s{+}1)$-tuple $(v_0,v_1,\ldots, v_s)$ of vertices with
$\{v_i,v_{i+1}\}\in E$ and $v_i\ne v_{i+2}$ for all possible $i$.
A $1$-arc is also called an arc.

Let $\Ga=(V,E)$ be a graph. A permutation $g$ on $V$ is called  an automorphism
of $\Ga$ if $\{u^g,v^g\}\in E$ for all $\{u,v\}\in E$. All automorphisms of $\Ga$
form a  subgroup of the symmetric group $\Sym(V)$, denoted by $\Aut\Ga$, which
is called the automorphism group of $\Ga$. The group $\Aut\Ga$ has a natural
action on   $E$, namely, $\{u,v\}^g=\{u^g,v^g\}$ for $\{u,v\}\in E$
and $g\in \Aut\Ga$. If this action is transitive, that is, for each
pair of  edges there
exists some $g\in \Aut\Ga$ mapping one   edge  to the other one, then
 $\Ga$ is called {\em edge-transitive}.
Similarly, we may define the {\em vertex-transitivity}, {\em arc-transitivity}
and {\em $s$-arc-transitivity} of $\Ga$.
The graph $\Ga$ is called {\em edge-primitive} if
$\Aut\Ga$ acts primitively on $E$, that is, $\Ga$ is edge-transitive and
the stabilizer $(\Aut\Ga)_{\{u,v\}}$ of some (and hence every) edge $\{u,v\}$
in $\Aut\Ga$ is a maximal subgroup.

Let $\Ga=(V,E)$ be an  edge-primitive graph of valency no less than $3$.
Then, as observed in \cite{Giu-Li-edge-prim}, $\Ga$ is also arc-transitive.
If $\Ga$ is $2$-arc-transitive then  Praeger's reduction theorems
\cite{Prag-o'Nan, Prag-o'Nan-bi} will be effective tools for us to
investigate the group-theoretic and graph-theoretic properties of $\Ga$.
However, $\Ga$ is not necessarily $2$-arc-transitive; for example,
by the Atlas \cite{Atlas},
the sporadic Rudvalis group $\Ru$
is the automorphism group of a rank $3$ graph, which is  edge-primitive
and of valency $2304$ but not $2$-arc-transitive.
Using  O'Nan-Scott Theorem for (quasi)primitive groups \cite{Prag-o'Nan},
Giudici and Li  \cite{Giu-Li-edge-prim} gave a reduction theorem on
the automorphism group of $\Ga$. They proved that,  as a primitive group on $E$,
only four of the eight O'Nan-Scott types for primitive groups may occur
for $\Aut\Ga$, say {\rm SD}, {\rm CD}, {\rm PA} and {\rm AS}.
They also considered the possible O'Nan-Scott types for $\Aut\Ga$ acting on $V$,
and presented constructions or examples to verify the existence of
corresponding graphs. Then what will happen if we assume further that
$\Ga$ is $2$-arc-transitive? The third author of this paper showed that either
$\Aut\Ga$ is almost simple or $\Ga$ is a complete bipartite graph if $\Ga$ is $2$-arc-transitive, see \cite{2-arc-edge-p}.
This stimulate our interest in classifying those edge-primitive graphs
which are  $2$-arc-transitive.

In this paper, we present a classification result stated as follows.
\begin{theorem}\label{main-result}
Let $\Ga=(V,E)$ be a  graph of valency $d\ge 6$, and let $G\le \Aut\Ga$ such that
 $G$ acts primitively on the edge set and transitively on the $2$-arc set of $\Ga$.
 Assume further that $G$ is almost simple and, for $\{u,v\}\in E$,  the edge-stabilizer $G_{\{u,v\}}$ is soluble.
 Then either  $\Ga$ is $(G,4)$-arc-transitive, or
 $G$, $G_{\{u,v\}}$, $G_v$ and $d$ are listed as in Table {\rm \ref{all-graphs}}.
\end{theorem}

\noindent{\bf Remark}.
 If $\Ga$ is  edge-primitive and either $4$-arc-transitive or of valency less
 than $6$, then the edge-stabilizers must be soluble. The reader may find a
 complete list of such graphs in \cite{5val, 4val, s-arc, 3val}.
 For each triple $(G, G_v, G_{\{u,v\}})$  listed   in
 Table {\rm \ref{all-graphs}}, the coset graph $\Cos(G, G_v, G_{\{u,v\}})$,
 see Section \ref{sect=pre} for the definition,
 is both $(G,2)$-arc-transitive and $G$-edge-primitive.
\qed

\begin{table}[ht]
\[
\begin{array}{|l|l|l|l|l|}
\hline
G & G_{\{u,v\}}& G_v& d& \mbox{Remark} \\ \hline

 \PSL_4(2).2 & 2^4{:}\S_4 & 2^3{:}\SL_3(2)&7&\\ \hline

 \PSL_5(2).2&[2^8]{:}\S_3^2.2& 2^6{:}(\S_3{\times} \SL_3(2))& 7& \\ \hline

 \F_4(2).2&[2^{22}]{:}\S_3^2.2 &[2^{20}]{.}(\S_3{\times} \SL_3(2))& 7& \\ \hline

  \PSL_4(3).2 & 3_{+}^{1{+}4}{:}(2\S_4{\times} 2)&3^3{:}\SL_3(3)&13&\\

 \PSL_4(3).2^2 & 3_{+}^{1{+}4}{:}(2\S_4{\times} \ZZ_2^2)&3^3{:}(\SL_3(3){\times}\ZZ_2)&13&\\ \hline

 \PSL_5(3).2& [3^8]{:}(2\S_4)^2.2&3^6{.}2\S_4.\SL_3(3)&13&\\ \hline \hline
\S_p & \ZZ_p{:}\ZZ_{p-1} & \PSL_2(p)& p{+}1& p\in \{7,11\}\\ \hline
\M_{11} & 3^2{:}\Q_8.2& \M_{10}& 10 &\K_{11}\\ \hline
\J_1& \ZZ_{11}{:}\ZZ_{10}& \PSL_2(11)& 12 &\\ \hline
\J_3.2& \ZZ_{19}{:}\ZZ_{18}& \PSL_2(19)& 20 &\\ \hline
\ON.2& \ZZ_{31}{:}\ZZ_{30}& \PSL_2(31)& 32 &\\ \hline
\B& \ZZ_{19}{:}\ZZ_{18}{\times} \ZZ_2& \PGL_2(19)& 20 &\\  \hline
\B &  \ZZ_{23}{:}\ZZ_{11}{\times} \ZZ_2& \PSL_2(23)& 24 &\\ \hline
 \M& \ZZ_{41}{:}\ZZ_{40} & \PSL_2(41)& 42 &\\ \hline
 \PSL_2(19)& \D_{20}& \PSL_2(5)& 6& \\ \hline
\A_6.2, \A_6.2^2 &  \ZZ_5{:}[4], \ZZ_{10}{:}\ZZ_{4}& \PSL_2(5),\PGL_2(5)& 6& \K_{6,6}, \, G\not\cong \S_6\\ \hline
 \PGL_2(11) & \D_{20}& \PSL_2(5)& 6&  \\ \hline
 \PSL_3(r) & 3^2{:}\Q_8& \PSL_2(9)& 10&  r \mbox{ is a prime with}  \\

  \PSL_3(r).2     & 3^2{:}\Q_8.2 &\PGL_2(9)& & r\equiv 4,16,31,34\,(\mod 45) \\ \hline

   \PSU_3(r) & 3^2{:}\Q_8& \PSL_2(9)& 10&  r \mbox{  is a prime with} \\

    \PSU_3(r).2&  3^2{:}\Q_8.2 & \PGL_2(9)& & r\equiv 11,14,29,41\,(\mod 45) \\ \hline
      \hline
\HS.2 & [5^3]{:}[2^5] & \PSU_3(5){:}2& 126&  \\ \hline
\Ru & [5^3]{:}[2^5] & \PSU_3(5){:}2& 126&  \\ \hline\hline

\M_{10}& \ZZ_8{:}\ZZ_2&3^2{:}\Q_8&9&\K_{10}\\ \hline
\PSL_3(3).2& \GL_2(3){:}2&3^2{:}\GL_2(3)&9& \\ \hline
\J_1& \ZZ_7{:}\ZZ_6& \ZZ_2^3{:}\ZZ_7{:}\ZZ_3& 8&\\ \hline

\PSL_2(p^f).[o]&\D_{2(p^f-1)\over (2,p-1)}.[o]& \ZZ_p^f{:}\ZZ_{p^f-1\over (2,p-1)}.[o]& p^f& \K_{p^f+1},\, o\div (2,p-1)f \\ \hline

\Sz(2^f).o&\D_{2(2^f-1)}.o&\ZZ_2^f{:}\ZZ_{2^f-1}.o& 2^f& f \mbox{ is odd, } o\div f\\ \hline
\end{array}
\]
{\caption{Graphs.} \label{all-graphs}}
\end{table}

\vskip 20pt

\section{Preliminaries}\label{sect=pre}

Let $G$ be a finite group and $H,K\le G$ with $|K:(H\cap K)|=2$ and
$\cap_{g\in G}H^g=1$. Let  $[G:H]=\{Hx\mid x\in G\}$, and define a
graph $\Cos(G,H,K)$ on $[G:H]$ such that
$\{Hx,Hy\}$ is an edge if and only if $yx^{-1}\in HKH\setminus H$.
The group $G$ can be viewed as a subgroup of $\Aut\Cos(G,H,K)$,
where $G$ acts on $[G:H]$ by right multiplication.
Then $\Cos(G,H,K)$ is $G$-arc-transitive and,
for  $x\in K\setminus H$, the edge $\{H,Hx\}$ has
stabilizer $K$ in $G$. Thus $\Cos(G,H,K)$
is $G$-edge-primitive if and only if $K$ is maximal in $G$.

Assume that $\Ga=(V,E)$ is a $G$-edge-primitive graph of valency $d\ge 3$.
Then $\Ga$ is $G$-arc-transitive by \cite[Lemma 3.4]{Giu-Li-edge-prim}.
Take an edge $\{u,v\}\in E$, let $H=G_v$ and $K=G_{\{u,v\}}$.
Then $K$ is maximal in $G$, and $H\cap K=G_{uv}$, which has index $2$ in $K$.
Noting that $\cap_{g\in G}H^g$ fixes   $V$ pointwise,  $\cap_{g\in G}H^g=1$.
Further,  $\,v^g\mapsto G_vg, \,\forall g\in G$ gives an isomorphism from $\Ga$ to
$\Cos(G,H,K)$. Then, by \cite[Theorem 2.1]{Fang99}, the following lemma holds.

\begin{lemma}\label{coset}
Let $\Ga=(V,E)$ be a connected graph of valency $d\ge 3$, and $G\le \Aut\Ga$.
Then $\Ga$ is both $(G,2)$-arc-transitive and $G$-edge-primitive if and only if
$\Ga\cong \Cos(G,H,K)$ for some subgroups $H$ and $K$ of $G$ satisfying
\begin{itemize}
\item[(1)] $|K:(H\cap K)|=2$, $\cap_{g\in G}H^g=1$ and $K$ is maximal in $G$;
\item[(2)] $H$ acts $2$-transitively on $[H:(H\cap K)]$ by right multiplication.
\end{itemize}
\end{lemma}

Let $\Ga=(V,E)$ be a connected graph of valency no less than $3$,
$\{u,v\}\in E$ and $G\le \Aut\Ga$. Assume that
$\Ga$ is $(G,s)$-arc-transitive for some $s\ge 1$, that is,
$G$ acts transitively on the $s$-arc set  of $\Ga$.
Then $G_v$ acts transitively on the neighborhood $\Ga(v)$ of $v$ in $\Ga$.
Let $\G_v^{\Ga(v)}$ be the transitive permutation group
induced by $G_v$ on  $\Ga(v)$, and
let $G_v^{[1]}$ be the kernel of $G_v$ acting on $\Ga(v)$.
Then $G_v^{\Ga(v)}\cong G_v/G_v^{[1]}$.
Set $G_{uv}^{[1]}=G_u^{[1]}\cap G_v^{[1]}$.
Then $G_v^{[1]}$ induces a normal subgroup of
$(G_{u}^{\Ga(u)})_v$ with the kernel $G_{uv}^{[1]}$.
Since $G$ is transitive on the arcs of $\Ga$,
there is some element in $G$ interchanging $u$ and $v$.
This implies that
\[|G_{\{u,v\}}{:}G_{uv}|=2 \mbox{ and } (G_v^{\Ga(v)})_u\cong (G_{u}^{\Ga(u)})_v.\]
Writing $G_v^{[1]}$ and $G_v$ in group extensions, the next lemma follows.
\begin{lemma}\label{extension}
\begin{enumerate}
\item[(1)] $G_v^{[1]}=G_{uv}^{[1]}.(G_v^{[1]})^{\Ga(u)}$, $(G_v^{[1]})^{\Ga(u)}\unlhd (G_{u}^{\Ga(u)})_v\cong  (G_v^{\Ga(v)})_u$.
\item[(2)] $G_{\{u,v\}}=G_{uv}.2$,  $G_{uv}=(G_{uv}^{[1]}.(G_v^{[1]})^{\Ga(u)}).(G_v^{\Ga(v)})_u$, $G_v=(G_{uv}^{[1]}.(G_v^{[1]})^{\Ga(u)}).G_v^{\Ga(v)}$.
    \item[(3)] If $G_{uv}^{[1]}=1$ then $G_{uv}\lesssim (G_v^{\Ga(v)})_u\times (G_{u}^{\Ga(u)})_v$.
\end{enumerate}
\end{lemma}

By \cite{Weiss81b}, $s\le 7$ , and if $s\ge 2$ then $G_{uv}^{[1]}$
is a $p$-group for some prime $p$, refer to \cite{Gardiner-73}.
Thus Lemma \ref{extension} yields a fact as follows.

\begin{corollary}\label{sl-estab}
Let $\Ga=(V,E)$ be a connected $(G,2)$-arc-transitive graph,
and $\{u,v\}\in E$. Then $G_{\{u,v\}}$ is soluble if and only if
$(G_v^{\Ga(v)})_u$  is soluble, and $G_v$ is soluble if and only if
$G_v^{\Ga(v)}$ is soluble.
\end{corollary}

Choose  $s$  maximal as possible, that is,    $\Ga$ is $(G,s)$-arc-transitive but not $(G,s{+}1)$-arc-transitive. In this case,   $\Ga$ is said to be
$(G,s)$-transitive. If further $G_{uv}^{[1]}\ne 1$,
then one can read out the vertex-stabilizer $G_v$ from
\cite{Gardiner74,Weiss81a} for $s\ge 4$ and
from \cite{Trofimov9} for $2\le s\le 3$.
In particular, we have the following result.
\begin{theorem}
\label{double-star}
Let $\Ga=(V,E)$ be a connected $(G,s)$-transitive graph of valency on less $3$,
and $\{u,v\}\in E$.
Assume that  $s\ge 2$.
\begin{enumerate}
\item[(1)] If $G_{uv}^{[1]}=1$ then $s=2$ or $3$.
\item[(2)] If $G_{uv}^{[1]}\ne 1$ then
 $G_{uv}^{[1]}$ is a  $p$-group for some prime $p$,
 $\PSL_n(q)\unlhd G_v^{\Ga(v)}$, $ |\Ga(v)|={q^n-1\over q-1}$ and $6\ne s\le 7$, where $n\ge 2$ and $q=p^f$ for some integer $f\ge 1$; moreover, either
\begin{enumerate}
\item[(2.1)] $n=2$ and $s\ge 4$; or
\item[(2.2)] $n\ge 3$, $s\le 3$ and   $\O_p(G_v)$  is given as in Table {\rm \ref{stab-struct}}, where $\O_p(G_v)$ is the maximal normal $p$-subgroup  of $G_v$.
\end{enumerate}
\end{enumerate}
\end{theorem}
\begin{table}[ht]
\[
\begin{array}{|l|l|l|l|l|l|}
\hline
\O_p(G_v)& G_{uv}^{[1]} & s & n & q& G_v  \\ \hline
\ZZ_p^{n(n-1)f}&\ZZ_p^{(n-1)^2f} & 3 &   & & \SL_{n-1}(q){\times} \SL_n(q)\unlhd G_v/\O_p(G_v) \\ \hline
 \ZZ_p^{nf} & \ZZ_p^f&  2 &   &   & a.\PSL_n(q)\unlhd G_v/\O_p(G_v)\mbox{ with } a\div q-1\\ \hline

\ZZ_p^{n(n-1)f\over 2}&\ZZ_p^{(n-1)(n-2)f\over 2} & 2 &   & & a.\PSL_n(q)\unlhd G_v/\O_p(G_v)\mbox{ with } a\div q-1\\ \hline
[q^{20}] & [q^{18}]& 3 & 3  & {\rm even} & \SL_2(q){\times} \SL_3(q)\unlhd G_v/\O_p(G_v) \\ \hline
\ZZ_3^6 & \ZZ_3^4 & 2&  3 & 3& \ZZ_3^6{:}\SL_3(3) \\ \hline
\ZZ_2^{n{+}1} &\ZZ_2^2& 2 &   & 2& \ZZ_2^{n{+}1}{:}\SL_n(2) \\ \hline
\ZZ_2^{11}, \ZZ_2^{14}& \ZZ_2^8, \ZZ_2^{11} & 2 &  4 & 2& \ZZ_2^{11}{:}\SL_4(2),\ZZ_2^{14}{:}\SL_4(2) \\ \hline
[2^{30}] & [2^{26}] &2 & 5  & 2& [2^{30}]{:}\SL_5(2) \\ \hline

\end{array}
\]
{\caption{}
\label{stab-struct}}
\end{table}

\begin{lemma}\label{O_r(G-alpha)}
Let $\Ga=(V,E)$ be a connected $(G,2)$-arc-transitive graph, and $\{u,v\}\in E$.
If $r$ is a prime divisor of $|\Ga(v)|$ then $\O_r(G_{v}^{[1]})=1=\O_r(G_{uv})$, and either
 $\O_r(G_v)=1$, or
$|\Ga(v)|=r^e$ and
 $\O_r(G_v)\cong \ZZ_r^e\cong \soc(G_v^{\Ga(v)})$ for some integer $e\ge 1$.
\end{lemma}
\proof
Since $\Ga$ is $(G,2)$-arc-transitive,
$G_v^{\Ga(v)}$ is a $2$-transitive group, and thus $G_{uv}$ is transitive on $\Ga(v)\setminus\{u\}$.
Since $\O_r(G_{uv}) \unlhd G_{uv}$, all $\O_r(G_{uv})$-orbits  on $\Ga(v)\setminus\{u\}$ have the same size.
Noting that  $r$ is coprime to $|\Ga(v)\setminus\{u\}|$, it follows that
$\O_r(G_{uv})\le G_v^{[1]}$. Since $G_v^{[1]} \unlhd G_{uv}$, we have $\O_r(G_v^{[1]})\le \O_r(G_{uv})$, and so $\O_r(G_v^{[1]})=\O_r(G_{uv})$.
Similarly, considering the action of $G_{uv}$ on $\Ga(u)\setminus\{v\}$, we get
$\O_r(G_u^{[1]})=\O_r(G_{uv})$. Then $\O_r(G_u^{[1]})=\O_r(G_{uv})=\O_r(G_v^{[1]})\le G_{uv}^{[1]}$.
By Theorem~\ref{double-star}, either  $G_{uv}^{[1]}=1$ or $G_{uv}^{[1]}$ is a $p$-group, where $p$ is a prime divisor of $|\Ga(v)|-1$. It follows that $\O_r(G_u^{[1]})=\O_r(G_{uv})=\O_r(G_v^{[1]})=1$.

Note that $\O_r(G_v)G_v^{[1]}/G_v^{[1]}\cong \O_r(G_v)/(\O_r(G_v)\cap G_v^{[1]})$.
Since $\O_r(G_v)\cap G_v^{[1]}\le \O_r(G_v^{[1]})=1$, we have $\O_r(G_v)\cong \O_r(G_v)G_v^{[1]}/G_v^{[1]}\unlhd G_v/G_v^{[1]}\cong \G_v^{\Ga(v)}$.
Thus   $\O_r(G_v)$ is isomorphic to a normal $r$-subgroup of $\G_v^{\Ga(v)}$.
This implies that either $\O_r(G_v)=1$, or $\G_v^{\Ga(v)}$ is an affine $2$-transitive group of degree $r^e$ for some $e$. Thus the lemma follows.
\qed

\vskip 5pt

Let $a\ge 2$ and $f\ge 1$ be  integers. A  prime divisor $r$ of $a^f-1$ is primitive if $r$ is not a divisor of $a^{e}-1$ for all $1\le e< f$. By  Zsigmondy's theorem \cite{Zs}, if  $f>1$ and $a^f-1$ has no primitive prime divisor then
$a^f=2^6$, or $f=2$ and $a=2^t-1$ for some prime $t$.
Assume that  $a^f-1$ has a primitive prime divisor $r$. Then $a$ has order $f$ modulo $r$. Thus $f$ is a divisor of $r-1$, and if $r$ is a divisor of $a^{f'}-1$ for some $f'\ge 1$ then $f$ is a divisor of $f'$.
Thus we have the following lemma.

\begin{lemma}\label{2-3} Let $a\ge 2$, $f\ge 1$ and $f'\ge 1$ be   integers. If
$a^f-1$ has a primitive prime divisor $r$ then $f$ is a divisor of $r-1$, and
$r$ is a divisor of  $a^{f'}-1$ if and only if $f$ is a divisor of $f'$.
If $f\ge 3$   then $a^f-1$ has a  prime divisor no less than $5$.
\end{lemma}

\vskip 5pt

We end this section with
a fact on finite primitive groups.

\begin{lemma}\label{sylow}
Assume that $G$ is a finite primitive group with a point-stabilizer $H$.
If $H$ has a normal Sylow subgroup $P\ne 1$, then $P$ is also a Sylow subgroup of $G$.
\end{lemma}
\proof
Assume that $P\ne 1$ is  a normal Sylow subgroup of $H$. Clearly, $P$ is not normal in $G$. Take a Sylow subgroup $Q$ of $G$ with $P\le Q$. Then $H\le \langle \N_Q(P),H\rangle\le \N_G(P)\ne G$. Since $H$ is maximal in $G$, we have
$H=\langle \N_Q(P),H\rangle$ and so $\N_Q(P)\le H$. It follows that $\N_Q(P)=P$, and hence $P=Q$. Then the lemma follows.
\qed

\vskip 20pt

\section{Some restrictions on  stabilizers}

In Sections \ref{sect-insoluble} and \ref{sect-soluble}, we shall prove Theorem \ref{main-result} using the   result given in \cite{s-arc} which classifies finite primitive groups with soluble point-stabilizers.
Let $\Ga=(V,E)$ be a   graph of valency $d\ge 6$,
$\{u,v\}\in E$ and $G\le \Aut\Ga$. Assume
that $G$ is almost simple, $G_{\{u,v\}}$ is soluble, $\Ga$ is $G$-edge-primitive and $(G,2)$-arc-transitive.
Clearly, each nontrivial normal subgroup of $G$ acts transitively on the edge set $E$.
Choose a minimal $X$ among the normal subgroups of $G$ which   act  primitively
on $E$. By the choice of $X$, we have $\soc(X)=\soc(G)$, $X_{\{u,v\}}=X\cap G_{\{u,v\}}$, $G=XG_{\{u,v\}}$ and
$G/X=XG_{\{u,v\}}/X\cong G_{\{u,v\}}/X_{\{u,v\}}$.
Then, considering  the  restrictions on both $X_{\{u,v\}}$ and $X_v$ caused by the
$2$-arc-transitivity of $\Ga$,
we may work out the pair $(X,X_{\{u,v\}})$ from \cite[Theorem 1.1]{s-arc}, and then determine the group $G$ and the graph $\Ga$.
Thus we make the following assumptions.

\begin{hypothesis}\label{hypo}
Let $\Ga=(V,E)$ be a $G$-edge-primitive graph of valency $d\ge 6$, and $\{u,v\}\in E$, where $G$ is  an almost simple group with socle $T$.  Assume that
\begin{itemize}
\item[(i)] $\Ga$ is $(G,2)$-arc-transitive, and the edge-stabilizer $G_{\{u,v\}}$ is  soluble;
\item[(ii)] $G$ has a normal subgroup $X$   such that $\soc(X)=T$, $X_{\{u,v\}}$ is maximal in $X$, and $(X,X_{\{u,v\}})$ is one of the pairs $(G_0,H_0)$  listed in
    \cite[Tables 14-20]{s-arc}.
 \end{itemize}
\end{hypothesis}

For the group $X$ in Hypothesis \ref{hypo}, we have $1\ne X_v^{\Ga(v)}\unlhd G_v^{\Ga(v)}$. Note that $G_v^{\Ga(v)}$ is $2$-transitive (on $\Ga(v)$).
Then $G_v^{\Ga(v)}$ is affine or almost simple, see \cite[Theorem 4.1B]{Dixon} for example. It follows that  $\soc(G_v^{\Ga(v)})=\soc(X_v^{\Ga(v)})$.

\subsection{}
Assume that $G_v$ is insoluble. Then $G_v^{\Ga(v)}$ is an almost simple $2$-transitive group (on $\Ga(v)$).  Recall that  $\soc(G_v^{\Ga(v)})=\soc(X_v^{\Ga(v)})$. Checking the point-stabilizers of almost simple $2$-transitive groups, since $(G_v^{\Ga(v)})_u$  is soluble, we conclude that either $X_v^{\Ga(v)}$ is $2$-transitive, or $G_v^{\Ga(v)}\cong \PSL_2(8){.}3$ and $d=28$.
For a complete list of finite $2$-transitive groups, the reader may refer to
\cite[Tables 7.3 and 7.4]{Cameron}.

\begin{lemma}\label{Xu-al-0}
Suppose  that Hypothesis {\rm\ref{hypo}} holds. If $d=28$ then  $G_v^{\Ga(v)}\not\cong \PSL_2(8){.}3$.
\end{lemma}
\proof
Suppose  that $G_v^{\Ga(v)}\cong \PSL_2(8){.}3$  and $d=28$.
Note that $X_{uv}^{[1]}\le G_{uv}^{[1]}=1$, see Theorem \ref{double-star}.
Thus $X_{uv}\lesssim (X_v^{\Ga(v)})_u{\times} (X_u^{\Ga(u)})_v$ by Lemma \ref{extension}.

Assume that $X_v^{\Ga(v)}\cong \PSL_2(8)$.
Then $(X_v^{\Ga(v)})_u\cong \D_{18}$, and $X_{uv}\cong \D_{18}$, $(\ZZ_3{\times}\ZZ_9){:}\ZZ_2$, $(\ZZ_9{\times}\ZZ_9){:}\ZZ_2$ or $\D_{18}{\times} \D_{18}$. In particular, the unique
Sylow $3$-subgroup of $X_{\{u,v\}}=X_{uv}.2$ is isomorphic to $\ZZ_m{\times} \ZZ_9$, where $m=1$, $3$ or $9$. Checking the primitive groups listed in \cite[Tables 14-20]{s-arc}, we know that only the pairs $(\PSL_2(q),\D_{2{q{\pm} 1\over (2,q-1)}})$ possibly meet our requirements on  $X_{\{u,v\}}$, yielding  $X_{\{u,v\}}\cong \D_{2{q{\pm} 1\over (2,q-1)}}$. Then $\D_{36}\cong X_{\{u,v\}}\cong \D_{2{q{\pm} 1\over (2,q-1)}}$. Calculation shows that $q=37$; however, $\PSL_2(37)$ has no subgroup which has  a quotient  $\PSL_2(8)$, a contradiction.

Now let $X_v^{\Ga(v)}=G_v^{\Ga(v)}\cong \PSL_2(8){.}3$. Then $(X_v^{\Ga(v)})_u\cong (X_u^{\Ga(u)})_v\cong \ZZ_9{:}\ZZ_6$ and $X_{uv}\lesssim \ZZ_9{:}\ZZ_6{\times} \ZZ_9{:}\ZZ_6$. In particular, a Sylow $2$-subgroup of $X_{\{u,v\}}=X_{uv}.2$
is not a cyclic group of order $8$, and the unique Sylow $3$-subgroup of
$X_{\{u,v\}}$  is nonabelian and contains elements of order $9$. Since $X_{\{u,v\}}=X_{uv}.2=X_v^{[1]}.(X_v^{\Ga(v)})_u.2$ and $X_v^{[1]}\cong (X_v^{[1]})^{\Ga(u)}\unlhd(X_u^{\Ga(u)})_v$,
we have $|X_{\{u,v\}}|=2^23^3$, $2^23^4$, $2^23^5$, $2^23^6$, $2^33^5$ or $2^33^6$.
Checking the Tables 14-20 given in \cite{s-arc}, we conclude that $X=\G_2(3).2$, and $X_{\{u,v\}}\cong [3^6]{:}\D_8$. In this case, $X_v^{[1]}\cong \ZZ_9{:}\ZZ_6$ and
$X_v\cong \ZZ_9{:}\ZZ_6.\PSL_2(8){.}3$; however, $X$ has no such   subgroup by the Atlas
\cite{Atlas},  a contradiction.
This completes the proof.
\qed

By Lemma \ref {Xu-al-0}, combining with Theorem \ref{double-star}, the next lemma
follows from checking the point-stabilizers of finite almost simple $2$-transitive groups.
\begin{lemma}\label{Xu-al-1}
Suppose  that Hypothesis {\rm\ref{hypo}} holds and $G_v^{\Ga(v)}$ is almost simple.  Then  one of the following holds{:}
\begin{enumerate}
\item[(s1)] $G_v^{\Ga(v)}=X_v^{\Ga(v)}=\PSL_3(2)$ or $\PSL_3(3)$, and $d=7$ or $13$, respectively;

    \item[(s2)] $\soc(X_v^{\Ga(v)})=\PSL_2(q)$ with $q>4$, and $d=q{+}1$;
    \item[(s3)] $G_{uv}^{[1]}=1$, $\soc(X_v^{\Ga(v)})=\PSU_3(q)$ with $q>2$, and $d=q^3{+}1$;
    \item[(s4)] $G_{uv}^{[1]}=1$,  $\soc(X_v^{\Ga(v)})=\Sz(q)$ with $q=2^{2n{+}1}>2$, and $d=q^2{+}1$;
    \item[(s5)]  $G_{uv}^{[1]}=1$, $\soc(X_v^{\Ga(v)}) =\Ree(q)$ with $q=3^{2n{+}1}>3$, and $d=q^3{+}1$.
\end{enumerate}
In particular, $\Ga$ is $(X,2)$-arc-transitive.
\end{lemma}

Recall that the Fitting subgroup $\Fit(H)$ of a finite group $H$ is the direct product of $\O_r(H)$, where $r$ runs over the set of prime divisors of $|H|$.
\begin{lemma}\label{Fitting}
Suppose  that Hypothesis {\rm\ref{hypo}} holds and one of Lemma {\rm \ref{Xu-al-1}~(s2)-(s5)} occurs. Let $q=p^f$ for some prime $p$.
Assume that $X_{u,v}^{[1]}=1$. Then $\Fit(X_{uv})=\O_p(X_{uv})$, and either
$\Fit(X_{uv})=\Fit(X_{\{u,v\}})$ or $\Fit(X_{\{u,v\}})=\Fit(X_{uv}).2$;
in particular,
$|\Fit(X_{\{u,v\}}):\O_p(X_{\{u,v\}})|\le 2$.
\end{lemma}
\proof
Let $r$ be a prime divisor of $|X_{uv}|$. Then $\O_r(X_{uv})$ is normal in $X_{uv}$. Since $\Ga$ is $(X,2)$-arc-transitive,
 $X_{uv}$ acts transitively on  $\Ga(v)\setminus \{u\}$. Thus all $\O_r(X_{uv})$-orbits (on $\Ga(v)\setminus \{u\}$) have equal size, which is a power of $r$ and a divisor of $|\Ga(v)\setminus \{u\}|$. Note that  $|\Ga(v)\setminus \{u\}|=d-1$, which
 is a power of $p$. It follows that either $r=p$ or $\O_r(X_{uv})=1$. Then $\Fit(X_{uv})=\O_p(X_{uv})$.

Note that $X_{uv}$ is normal in $X_{\{u,v\}}$ as $|X_{\{u,v\}}:X_{uv}|=2$.
Since $\O_p(X_{uv})$ is a characteristic subgroup of $X_{uv}$, it follows that  $\O_p(X_{uv})$ is normal in $X_{\{u,v\}}$, and so $\O_p(X_{uv})\le \O_p(X_{\{u,v\}})\le \Fit(X_{\{u,v\}})$.
For each odd prime divisor $r$ of $|X_{\{u,v\}}|$, since $|X_{\{u,v\}}:X_{uv}|=2$,
 we have $\O_r(X_{\{u,v\}})\le X_{uv}$, and so $\O_r(X_{\{u,v\}})=\O_r(X_{uv})$.
 It follows that $\Fit(X_{\{u,v\}})=\Fit(X_{uv})\O_2(X_{\{u,v\}})=\O_p(X_{uv})\O_2(X_{\{u,v\}})$.
 In particular, $\O_p(X_{uv})=\O_p(X_{\{u,v\}})$ if $p\ne 2$.

 It is easily shown that $X_{uv}\cap \O_2(X_{\{u,v\}})=\O_2(X_{uv})$.
 If $X_{uv}\ge \O_2(X_{\{u,v\}})$ then $p=2$, $\Fit(X_{\{u,v\}})=\O_2(X_{\{u,v\}})=\Fit(X_{uv})$, and the lemma holds.
 Assume that $\O_2(X_{\{u,v\}})\not\le X_{uv}$.
Since $|X_{\{u,v\}}:X_{uv}|=2$, we have $X_{\{u,v\}}=X_{uv}\O_2(X_{\{u,v\}})$.
Then $2|X_{uv}|=|X_{\{u,v\}}|=|X_{uv}||\O_2(X_{\{u,v\}}):(X_{uv}\cap \O_2(X_{\{u,v\}}))|=|X_{uv}||\O_2(X_{\{u,v\}}):\O_2(X_{uv})|$, yielding
$|\O_2(X_{\{u,v\}}):\O_2(X_{uv})|=2$. If $p=2$ then $\Fit(X_{\{u,v\}})=\O_2(X_{\{u,v\}})$ and $\Fit(X_{uv})=\O_2(X_{uv})$,  the lemma follows. If $p\ne 2$ then $\O_2(X_{uv})=1$,
$|\O_2(X_{\{u,v\}})|=2$, and so $\Fit(X_{\{u,v\}})=\O_p(X_{uv})\times \ZZ_2$.
This completes the proof.
\qed

\subsection{} Assume that  Hypothesis {\rm\ref{hypo}} holds and $G_v$ is soluble.
Then $G_v^{\Ga(v)}$ is an affine $2$-transitive group. Let $\soc(G_v^{\Ga(v)})=\ZZ_p^f$. Then $d=p^f$.
Recalling that $d\ge 6$, we have $G_{uv}^{[1]}=1$ by Theorem \ref{double-star}, and so $G_{uv}\lesssim(G_v^{\Ga(v)})_u{\times} (G_u^{\Ga(u)})_v$.
If $G_{uv}$ is abelian then $\Ga$ is known by   \cite{2-arc-edge-p}.
Thus we assume further that $G_{uv}$ is not abelian.
Then $(G_v^{\Ga(v)})_u$ is nonabelian,
and so $(G_v^{\Ga(v)})_u\not\le \GL_1(p^f)$; in particular, $f>1$. Since $(G_v^{\Ga(v)})_u$ is soluble, by \cite[Table 7.3]{Cameron},
we have the following lemma.
\begin{lemma}\label{stab-soluble-1}
Suppose that  Hypothesis {\rm\ref{hypo}} holds, $G_v$ is soluble and $G_{uv}$ is not abelian. Let $\soc(G_v^{\Ga(v)})=\ZZ_p^f$, where $p$ is a prime. Then $f>1$, and one of the following holds.
\begin{enumerate}
\item[(a1)]  $f=2$, and either $\SL_2(3)\unlhd (G_v^{\Ga(v)})_u\le \GL_2(p)$ and $p\in \{3,5,7,11,23\}$, or $p=3$ and $(G_v^{\Ga(v)})_u=\Q_8$;
\item[(a2)] $2_{+}^{1{+}4}{:}\ZZ_5\le (G_v^{\Ga(v)})_u\le 2_{+}^{1{+}4}.(\ZZ_5{:}\ZZ_4)< 2_{+}^{1{+}4}.\S_5$,
and $p^f=3^4$;
\item[(a3)] $(G_v^{\Ga(v)})_u\not\le \GL_1(p^f)$, $(G_v^{\Ga(v)})_u\le\GammaL_1(p^f)$ and $|(G_v^{\Ga(v)})_u|$ is divisible by $p^f-1$.
\end{enumerate}
\end{lemma}

Consider the case (a3) in Lemma \ref{stab-soluble-1}.
Write \[\GammaL_1(p^f)=\langle \tau,\sigma\mid  \tau^{p^f-1}=1=\sigma^f, \sigma^{-1}\tau\sigma=\tau^p\rangle.\]
Let $\langle \tau\rangle\cap (G_v^{\Ga(v)})_u=\langle \tau^m\rangle$, where $m\div (p^f-1)$.
Then \[(G_v^{\Ga(v)})_u/\langle \tau^m\rangle\cong \langle \tau\rangle  (G_v^{\Ga(v)})_u/\langle \tau\rangle\lesssim\langle\sigma\rangle.\]
Set $(G_v^{\Ga(v)})_u/\langle \tau^m\rangle\cong \langle\sigma^e\rangle$ for some divisor  $e$   of $f$. Then \[(G_v^{\Ga(v)})_u\cong \ZZ_{p^f-1\over m}.\ZZ_{f\over e}.\]
Choose $\tau^l\sigma^k\in (G_v^{\Ga(v)})_u$ such that $(G_v^{\Ga(v)})_u=\langle \tau^m\rangle\langle \tau^l\sigma^k\rangle$. Then $(\tau^l\sigma^k)^{f\over e}\in \langle \tau^m\rangle$ but $(\tau^l\sigma^k)^j\not\in \langle \tau^m\rangle$ for $1\le j< {f\over e}$. It follows that $\sigma^k$ has order $f\over e$.
Then $\sigma^k=\sigma^{ie}$ for some $i$ with $(i, {f\over e})=1$, and then $(\sigma^k)^{i'}=\sigma^e$ for some $i'$. Thus, replacing $\tau^l\sigma^k$ by   a power of it if necessary, we may let $k=e$. Then \[(G_v^{\Ga(v)})_u=\langle \tau^m\rangle\langle \tau^l\sigma^e\rangle.\]
Further, $(G_v^{\Ga(v)})_u=\langle \tau^m\rangle\langle (\tau^m)^i\tau^l\sigma^e\rangle$ for an arbitrary integer $i$, thus we may assume further $0\le l<m$. By \cite[Proposition 15.3]{Foulser}, letting ${\rm \pi}(n)$ be the set of prime divisors of a positive integer $n$, we have
\begin{enumerate}
\item[($\divideontimes$)] ${\rm \pi}(m)\subseteq {\rm \pi}(p^e-1)$, $me\div f$ and $(m,l)=1$; in particular, $m=1$ if $l=0$.
\end{enumerate}

Suppose that $X_{uv}$ is nonabelian. (The case where $X_{uv}$ is  abelian is left in Section \ref{sect-soluble}.)
Since $X_{uv}^{[1]}\le \G_{uv}^{[1]}=1$, we have
 \[X_v^{[1]}\unlhd (X_u^{\Ga(u)})_v\cong (X_v^{\Ga(v)})_u,\,\,
 X_{uv}\lesssim (X_u^{\Ga(u)})_v{\times}  (X_v^{\Ga(v)})_u.\]
This yields that $(X_v^{\Ga(v)})_u$ is nonabelian. Then a limitation on $\pi(|X_{uv}|)$ is given as follows.

\begin{lemma}\label{(primitive-di)}
Assume that Lemma {\rm \ref{stab-soluble-1}~(a3)} holds and $X_{uv}$ is nonabelian. Then  $(X_v^{\Ga(v)})_u\cong \ZZ_{m'}.\ZZ_{f\over e'}$, where $m'$ and $e'$ satisfy
\begin{enumerate}
\item[(i)] $\ZZ_{m'}\cong (X_v^{\Ga(v)})_u\cap \langle\tau^m\rangle$, $mm'\div p^f-1$, $e\div e'\div f$; and
    \item[(ii)]  $m'>1$, $e'<f$, ${\rm \pi}(p^f-1)\setminus {\rm \pi}(p^{e'}-1)\subseteq{\rm \pi}(m')\subseteq \pi(|X_{uv}|)$.
\end{enumerate}
\end{lemma}
\proof
Recall that $(X_v^{\Ga(v)})_u\unlhd (G_v^{\Ga(v)})_u=\langle \tau^m\rangle\langle \tau^l\sigma^e\rangle\cong \ZZ_{p^f-1\over m}.\ZZ_{f\over e}$. Then \[(X_v^{\Ga(v)})_u/((X_v^{\Ga(v)})_u\cap \langle \tau^m\rangle)\cong (X_v^{\Ga(v)})_u\langle \tau^m\rangle/ \langle \tau^m\rangle\lesssim \ZZ_{f\over e},\]
yielding $(X_v^{\Ga(v)})_u\cong \ZZ_{m'}.\ZZ_{f\over e'}$ with $m'$ and $e'$ satisfying (i). Since  $X_{uv}$ is nonabelian, $(X_v^{\Ga(v)})_u$ is nonabelian, and so $m'>1$ and $e'<f$.

By the above $(\divideontimes)$, each $r\in {\rm \pi}(p^f-1)\setminus {\rm \pi}(p^{e'}-1)$ is a divisor of $|\langle \tau^m\rangle|={p^f-1\over m}$. Let $R$ be the unique subgroup of order $r$ of $\langle \tau^m\rangle$. Then, noting that $R$ is normal in $(G_v^{\Ga(v)})_u$, either $R\le (X_v^{\Ga(v)})_u$ or $R (X_v^{\Ga(v)})_u=R{\times}(X_v^{\Ga(v)})_u$. Suppose that the latter case occurs. Since $e'<f$, we may let $\tau^n\sigma^{e'}\in (X_v^{\Ga(v)})_u\setminus \langle \tau^m\rangle$.
Then $\sigma^{e'}$ centralizes $R$.
Thus $x^{p^{e'}}=x$ for $x\in R$, yielding $r\div (p^{e'}-1)$,  a contradiction.
Then $R\le (X_v^{\Ga(v)})_u\cap \langle\tau^m\rangle\cong \ZZ_{m'}$.
Noting that $m'$ is a divisor of $|X_{uv}|$,   the result follows.
\qed

 \vskip 20pt

\section{Graphs with insoluble vertex-stabilizers}\label{sect-insoluble}

In this and next sections, we   prove  Theorem \ref{main-result}. Thus, we let $G$, $T$, $X$ and $\Ga=(V,E)$ be as in  Hypothesis \ref{hypo}.
Our task is to determine which pair $(G_0,H_0)$ listed in \cite[Tables 14-20]{s-arc} is a possible candidate for $(X,X_{\{u,v\}})$, and determine whether or not the
resulting triple $(G,G_v,G_{\{u,v\}})$ meets the conditions (1) and (2) in Lemma \ref{coset}.

In this section, we   deal with the case where $G_v$ is insoluble, that is,
$X_v$ is described as in Lemma \ref{Xu-al-1}. First, by the following lemma,
Lemma \ref{Xu-al-1} (s4) and (s5) are excluded.

\begin{lemma}\label{s4-s5}
Lemma {\rm \ref{Xu-al-1}} {\rm (s4)} and {\rm (s5)}  do not occur.
\end{lemma}
\proof
Suppose that Lemma \ref{Xu-al-1} (s4) or (s5)  holds.
Then  $X_{uv}^{[1]}=1$ by Theorem \ref{double-star}. Thus
$X_v=X_v^{[1]}.X_v^{\Ga(v)}$, $X_v^{[1]}\cong (X_v^{[1]})^{\Ga(u)}\unlhd (X_u^{\Ga(u)})_v\cong (X_v^{\Ga(v)})_u$,
 and
$X_{uv}\lesssim (X_v^{\Ga(v)})_u{\times} (X_u^{\Ga(u)})_v$. Set $q=p^f$ with $p$ a prime. Then the pair $(X_v^{\Ga(v)}, (X_v^{\Ga(v)})_u)$ is given as follows{:}
\[
\begin{array}{l|l|l}
X_v^{\Ga(v)}  & (X_v^{\Ga(v)})_u  & \\ \hline
\Sz(q).e & p^{f{+}f}{:}(q-1).e&  e\mbox{ a divisor of } f,\,p=2,\, \mbox{odd } f>1 \\ \hline
\Ree(q).e &  p^{f{+}2f}{:}(q-1).e &e\mbox{ a divisor of } f,\,p=3,\, \mbox{odd } f>1
\end{array}
\]
In particular,  $\O_p(X_{\{u,v\}})$ is not abelian.

We next show that  none of the pairs $(G_0,H_0)$  in \cite[Tables 14-20]{s-arc}  gives  a desired pair $(X,X_{\{u,v\}})$.
Since $\O_p(X_{\{u,v\}})$ is nonabelian, those pairs  $(G_0,H_0)$ with $\O_p(H_0)$ abelian are not in our consideration. In particular,  $\soc(X)$ is not isomorphic to an alternating group. Also, noting that $X_{\{u,v\}}$ has a subgroup of index $2$, those $H_0$ having no  subgroup of index $2$ are excluded.

{\bf Case 1}.
Suppose that $\soc(X_v^{\Ga(v)})=\Ree(q)$. Then $p=3$,
$\O_3(X_{\{u,v\}})$ is nonabelian and of order $3^{3f}$, $3^{4f}$, $3^{5f}$ or $3^{6f}$, $|X_{\{u,v\}}|$ is a divisor of $2(q-1)^2f^2$ and divisible by $2(q-1)$. Checking the orders of those $H_0$ given in \cite[Tables 15]{s-arc}, we conclude  that $\soc(X)$ is not a sporadic simple group.

Suppose that $\soc(X)$ is a simple exceptional group of Lie type.
By \cite[Table 20]{s-arc}, we conclude that $(X,X_{\{u,v\}})$ is one of  $(\G_2(3^t).\ZZ_{2^{l+1}}, [3^{6t}]{:}\ZZ_{3^t-1}^2.\ZZ_{2^{l+1}})$ and $(\Ree(3^t),[3^{3t}]{:}\ZZ_{3^t-1})$, where $2^l$ is the $2$-part  of $t$. If $f>t$ then, by Zsigmondy's theorem, $q-1$ has a prime divisor which is
not a divisor of $3^t-1$, which contradicts that $|X_{\{u,v\}}|$ is   divisible by $2(q-1)$. Thus $f=t$, and then $X=\G_2(q).\ZZ_{2^{l+1}}$ and $X_{\{u,v\}}\cong [q^6]{:}\ZZ_{q-1}^2.\ZZ_{2^{l+1}}$.
This implies that $X_v^{[1]}\ne 1$, in fact, $|\O_3(X_v^{[1]})|=q^3$.
Thus $\O_3(X_v)\ne 1$ and $X_v$ has a quotient $\Ree(q).e$.
Checking the maximal subgroups of $\G_2(q).\ZZ_{2^{l+1}}$, refer to \cite[Theorems A and B]{Kleidman-2}, we conclude that $\G_2(q).\ZZ_{2^{l+1}}$ has no  maximal subgroup containing such $X_v$ as a subgroup, a contradiction.

Suppose that $\soc(X)$  is a  simple  classical group
over a finite field of order $r^t$, where $r$ is a prime.
Since $f>1$ is odd, $3^f-1$ has an odd prime divisor, and so $X_{\{u,v\}}$ is not a
$\{2,3\}$-group as $|X_{\{u,v\}}|$ is divisible by $3^f-1$.
Recall that $\O_3(X_{\{u,v\}})$ is nonabelian and of order $3^{3f}$, $3^{4f}$, $3^{5f}$ or $3^{6f}$. Checking  the groups $H_0$ given in \cite[Table 16-19]{s-arc},
we conclude that $\soc(X)=\PSL_n(r^t)$ or $\PSU_n(r^t)$, where $n\in \{3,4\}$.
Take a maximal subgroup $M$ of $X$ such that $X_v\le M$. Then $M$ has a simple section (i.e., a quotient of some subgroup) $\Ree(q)$. Recall that $q>3$. Checking Tables 8.3-8.6 and 8.8-8.11 given in \cite{Low}, we conclude that none of $\PSL_3(r^t)$, $\PSL_4(r^t)$, $\PSU_3(r^t)$ and $\PSU_4(r^t)$ has such maximal subgroups, a contradiction.

{\bf Case 2}.
Suppose that $\soc(X_v^{\Ga(v)})=\Sz(q)$. Then $q=2^f$,
$|\O_2(X_{\{u,v\}})|=2^{2f}a$, $2^{3f}a$ or $2^{4f}a$, where $f>1$ is odd, and $a=1$ or $2$.
Noting that $X_{\{u,v\}}$ has order divisible by $2^f-1$, by Lemma \ref{2-3}, we conclude that $X_{\{u,v\}}$ is not a $\{2,3\}$-group. Since $X_{\{u,v\}}$
is nonabelian, it follows from \cite[Table 15-20]{s-arc} that  either
$(X,X_{\{u,v\}})$ is one of $({}^2\F_4(2)', [2^9]{:}5{:}4)$, $(\PSp_4(2^t).\ZZ_{2^{l+1}}, [2^{4t}]{:}\ZZ_{2^t-1}^2.\ZZ_{2^{l+1}})$ and $(\Sz(2^t),[2^{2t}]{:}\ZZ_{2^t-1})$, or
  $\soc(X)$ is one of $\PSL_n(r^t)$ and $\PSU_n(r^t)$,  where $n\in \{3,4\}$, $2^l$ is the $2$-part of $t$, and $r$ is odd if $n=4$.
The first pair  leads to $q=2^3$, and so $|X_{\{u,v\}}|$ is divisible by $7$, a contradiction. Checking   the maximal subgroups   of $\soc(X)$ (refer to \cite[Tables 8.3-8.6, 8.8-8.14]{Low}), the groups   $\PSL_3(r^t)$,  $\PSU_3(r^t)$, $\PSL_4(r^t)$ and  $\PSU_4(r^t)$ are excluded as they have no maximal subgroup with a simple section $\Sz(q)$. Thus  $(X,X_{\{u,v\}})=(\PSp_4(2^t).\ZZ_{2^{l+1}}, [2^{4t}]{:}\ZZ_{2^t-1}^2.\ZZ_{2^{l+1}})$ or $(\Sz(2^t),[2^{2t}]{:}\ZZ_{2^t-1})$.

Since $f>1$ is odd, $2^f-1$ has a primitive prime divisor, say  $s$.
Recalling that $|X_{\{u,v\}}|$ is  divisible by $2^f-1$, it follows that
$s$ is a divisor of $2^t-1$, and so $t\ge f$. Then $t=f$. It follows that $X=\PSp_4(q).\ZZ_{2^{l+1}}$, and
$X_v^{[1]}\cong [q^2]{:}\ZZ_{q-1}$. However, by  \cite[Table 8.14]{Low}, $\PSp_4(q).\ZZ_{2^{l+1}}$ has no maximal subgroup containing $[q^2]{:}\ZZ_{q-1}.\Sz(q)$, a contradiction.
\qed

\begin{lemma}\label{d=7-13}
Assume that {\rm (s1)} of Lemma {\rm \ref{Xu-al-1}} occurs. Then $G$, $X$, $X_{\{u,v\}}$ and $X_v$ are listed as in Table {\rm \ref{small-gh}}.
\begin{table}[ht]
\[
\begin{array}{|l|l|l|l|l|l|}
\hline
G& X& X_{\{u,v\}}& X_v & s & d  \\ \hline
X &\PSL_4(2).2, \S_8 & 2^4{:}\S_4 & 2^3{:}\SL_3(2)&2&7\\ \hline
X& \PSL_5(2).2&[2^8]{:}\S_3^2.2& 2^6{:}(\S_3{\times} \SL_3(2))& 3& 7\\ \hline
X&\F_4(2).2&[2^{22}]{:}\S_3^2.2 &[2^{20}]{.}(\S_3{\times} \SL_3(2))& 3& 7\\ \hline
X,X.2& \PSL_4(3).2 & 3_{+}^{1{+}4}{:}(2\S_4{\times} 2)&3^3{:}\SL_3(3)&2&13\\ \hline
X& \PSL_5(3).2& [3^8]{:}(2\S_4)^2.2&3^6{.}2\S_4.\SL_3(3)&3&13\\ \hline
\end{array}
\]
{\caption{}
\label{small-gh}}
\end{table}
\end{lemma}
\proof
Assume first that $X_{uv}^{[1]}=1$. Then,
$X_v=X_v^{[1]}.X_v^{\Ga(v)}$, $X_v^{[1]}\cong (X_v^{[1]})^{\Ga(u)}\unlhd (X_u^{\Ga(u)})_v\cong (X_v^{\Ga(v)})_u$,
 and
$X_{uv}\lesssim (X_v^{\Ga(v)})_u{\times} (X_u^{\Ga(u)})_v$.

Suppose that $X_v^{\Ga(v)}=\PSL_3(2)$. Then $(X_v^{\Ga(v)})_u\cong \S_4$, and thus $X_v^{[1]}$ and $X_{\{u,v\}}$ are given as follows{:}
\[
\begin{array}{l|l|l|l|l}
X_v^{[1]}& 1 & 2^2 &\A_4&\S_4\\ \hline
X_{\{u,v\}}& 2^2{:}\S_3.2& 2^4.\S_3.2&2^4{:}3^2.[4]& 2^4{:}\S_3^2.2\\
\end{array}
\]
In particular, $2^2\le |\O_2(X_{\{u,v\}})|\le 2^5$. Check all possible pairs $(X,X_{\{u,v\}})$ in \cite[Tables 14-20]{s-arc}. Noting that
$\A_8\cong \PSL_4(2)$ and $\PSU_4(2)\cong \PSp_4(3)$, we conclude  that
$X\cong \A_8$, $X_{\{u,v\}}\cong 2^4{:}\S_3^2$ and $X_v^{[1]}\cong \A_4$; or
$X= \M_{12}$ with $X_{\{u,v\}}\cong 2_{+}^{1{+}4}{:}\S_3$; or
$X\cong \PSU_4(2)$ with $X_{\{u,v\}}\cong 2\A_4^2.2$. The group $\A_8$ is excluded as it has no   subgroup of the form of
$X_v^{[1]}.\PSL_3(2)$. The groups $\M_{12}$ and $\PSU_4(2)$
are excluded as their orders are not divisible by $d=7$.

Suppose that $X_v^{\Ga(v)}=\PSL_3(3)$. Then $(X_v^{\Ga(v)})_u\cong 3^2{:}2\S_4$.
Thus $X_v^{[1]}$ and $X_{\{u,v\}}$ are given as follows{:}
\[
\begin{array}{c|l|l|l|l|l|l}
X_v^{[1]}& 1 & 3^2 &3^2{:}2&3^2\Q_8& 3^2{:}2\A_4& 3^2{:}2\S_4\\ \hline
X_{\{u,v\}}& 3^2{:}2\S_4{.}2& 3^4{:}2\S_4{.}2& 3^4{:}([4].\S_4){.}2&3^4{:}\Q_8^2{.}\S_3{.}2& 3^4{:}(2\A_4)^2{.}[4] &3^4{:}(2\S_4)^2{.}2\\
\end{array}
\]
Note that $\O_3(X_{\{u,v\}})\cong 3^2$ or $3^4$. Checking the possible pairs $(X,X_{\{u,v\}})$, we have $X_{\{u,v\}}\cong 3^4{:}2^3{.}\S_4$ and
$X= \A_{12}$ or $\POmega^{+}_8(2)$; in this case, $d=13$ is not a divisor of
$|X|$, a contradiction.

Now let $X_{uv}^{[1]}$ be a nontrivial $p$-group.  Then, by Theorem \ref{double-star},
 $X_v$ and $X_{\{u,v\}}$ are given as follows{:}
\[
\begin{array}{c|c|l|l|l}
X_v & X_{\{u,v\}}& s &   d  & p\\ \hline
2^6{.}(\S_3{\times}\SL_3(2)) & [2^{8}]{.}\S_3^2.2& 3 & 7 & 2\\ \hline
[2^{20}]{.}(\S_3{\times}\SL_3(2))& [2^{22}]{.}\S_3^2.2 & 3 & 7&  2\\ \hline
 2^3{.}\SL_3(2) &[2^5].\S_3.2&   2 &    7& 2\\ \hline
 2^4{:}\SL_3(2)  &[2^6].\S_3.2&   2 &    7& 2\\ \hline
3^6{.}(2\A_4{\times}\SL_3(3))&[3^8]{.}(2\A_4{\times} 2\S_4).2&  3 & 13  & 3\\
3^6{.}(2\S_4{\times}\SL_3(3))&[3^8]{.}(2\S_4)^2.2 & 3 & 13  & 3\\ \hline
 3^3{.}\SL_3(3)&[3^5]{.}2\S_4.2 &   2 &    13 & 3\\
 3^3{.}(2{\times} \SL_3(3))&[3^5]{.}(2{\times} 2\S_4).2 &   2 &    13 & 3\\ \hline
 3^6{:}\SL_3(3) &[3^8]{.}2\S_4.2&   2 &    13 &3 \\
\end{array}
\]

Suppose that $p=2$. Then $|X_{\{u,v\}}|$ is divisible by $9$  if and only if $|\O_2(X_{\{u,v\}})|\ge 8$, and $\O_2(X_{\{u,v\}})$ contains no elements of order $8$ unless $|\O_2(X_{\{u,v\}})|\ge 2^{22}$. Check the pairs $(G_0,H_0)$ given in \cite[Tables 14-20]{s-arc} by estimating $|H_0|$ and $|\O_2(H_0)|$. We conclude that
one of the following holds{:}
\begin{enumerate}
\item[(i)] $X=\PSL_4(2).2\cong \S_8$ and $X_{\{u,v\}}= 2^4{:}\S_4$;
\item[(ii)] $X= \PSL_5(2).2$ and $X_{\{u,v\}}= [2^8]{.}\S_3^2.2$;
\item[(iii)] $X= \F_4(2).2$ and $X_{\{u,v\}}= [2^{22}]{.}\S_3^2.2$;
    \item[(iv)] $\soc(X)=\PSL_3(4)$ and $|\O_2(X_{\{u,v\}})|=2^6$;
     \item[(v)] $X= \PSU_4(3).2_3$ and $|\O_2(X_{\{u,v\}})|=2^7$;
     \item[(vi)] $X= \He.2$ and $X_{\{u,v\}}= [2^8]{:}\S_3^2.2$.
\end{enumerate}
Case (iv) yields that $X_v\cong 2^3{:}\SL_3(2)$ or $2^4{:}\SL_3(2)$; however, $X$
has no such  subgroup by the Atlas \cite{Atlas}. Similarly, cases (v) and (vi) are excluded.
For (i),
$G=X$  and $\Ga$ is (isomorphic to) the
point-plane incidence graph of the projective geometry $\PG(3,2)$.
For (ii),
 $G=X$ and $\Ga$ is (isomorphic to) the
line-plane incidence graph of the projective geometry $\PG(4,2)$.
If (iii) holds then $G=X$ and $\Ga$ is the line-plane incidence graph of the metasymplectic space associated with $\F_4(2)$, see \cite{Weiss-pro}.

Now let $p=3$. Then $|\O_3(X_{\{u,v\}})|=3^5$ or $3^8$, and $X_{\{u,v\}}$ has no normal Sylow subgroups. Checking all possible pairs $(X,X_{\{u,v\}})$ in
\cite[Tables 14-20]{s-arc},
we know that $(X,X_{\{u,v\}})$ is one of the following pairs:
\[(\F_4(8).2, 9^4.(2_{+}^{1{+}4}{:}\S_3^2).2),\,(\PSL_5(3).2, [3^8]{:}(2\S_4)^2.2), \,(\PSL_4(3).2, 3_{+}^{1{+}4}{:}(2{\times} 2\S_4)).\]
Note that $\O_3(X_v)\le \O_3(X_{\{u,v\}})$.
 Then, for the first pair,
 $\O_3(X_{\{u,v\}})\cong \ZZ_9^4$ has no   subgroup isomorphic to $\ZZ_3^6$, which is impossible. For the second pair,  $G=X$  and  $\Ga$ is (isomorphic to) the
line-plane incidence graph of the projective geometry $\PG(4,3)$.
The last pair implies that $X\not\cong \PGL_4(3)$, $G=X$ or $X.2$, and $\Ga$ is (isomorphic to) the
line-plane incidence graph of the projective geometry $\PG(3,3)$.
\qed

\begin{lemma}\label{PSL_2(q)}
Assume that   Lemma {\rm \ref{Xu-al-1}~(s2)} holds.
Then $d=q+1$, and either $\Ga$ is $(X,4)$-arc-transitive, or $G$, $X$, $X_{\{u,v\}}$ and $X_v$ are listed in Table {\rm \ref{PSL2q-graph}}.
\end{lemma}
\proof
Let $X_v^{\Ga(v)}= \PSL_2(q).[o]$, and $q=p^f> 4$, where $p$ is a prime and $o\div (2,q-1)f$.
Note that $\Ga$ is $(X,2)$-arc-transitive, see Lemma \ref{Xu-al-1}.
By Theorem \ref{double-star}, if $X_{uv}^{[1]}\ne 1$ then  $\Ga$ is $(X,4)$-arc-transitive. Thus we assume next that $X_{uv}^{[1]}=1$, and then Lemma \ref{Fitting} works.

Note that
$X_v=X_v^{[1]}.X_v^{\Ga(v)}$, $X_v^{[1]}\cong (X_v^{[1]})^{\Ga(u)}\unlhd (X_u^{\Ga(u)})_v\cong (X_v^{\Ga(v)})_u= p^f{:}{q-1\over (2,q-1)}{.}[o]$,
 and
$X_{uv}\lesssim (X_v^{\Ga(v)})_u{\times} (X_u^{\Ga(u)})_v$.  We have
$\O_p(X_{\{u,v\}})=\ZZ_p^{if}.a$, where $i\in \{1,2\}$ and $a$ is a divisor of $(2,p)$.
It is easily shown that $i=2$ if and only if $\O_p(X_v^{[1]})=\ZZ_p^f$.
Combining with Lemma \ref{Fitting},  we need only consider those pairs $(G_0,H_0)$ in
\cite[Tables 14-20]{s-arc} which satisfy
\begin{enumerate}
\item[(a)]  $\O_p(H_0)=\ZZ_p^{if}.a$, where $i\in \{1,2\}$ and $a$ is a divisor of $(2,p)$;
$|\Fit(H_0):\O_p(H_0)|\le 2$; $G_0$ has a subgroup, say $M_0$,  such that $|M_0:(M_0\cap H_0)|=q+1$,
$|H_0:(M_0\cap H_0)|=2$, and $M_0$ has a simple section $\PSL_2(q)$;

\item[(b)] $|H_0{:}\O_p(H_0)|$ is a divisor of $2(q-1)^2f^2$ and   divisible by $q-1$; if $i=1$ then $|H_0{:}\O_p(H_0)|$ is a divisor of $2(q-1)f$.
\end{enumerate}

{\bf Case 1}. Assume that $\soc(X)$ is an alternating group.
Using \cite[Table 14]{s-arc},
we have $G=X=\S_p$ and $X_{\{u,v\}}\cong \ZZ_p{:}\ZZ_{p-1}$, where $p\in \{7,11,17,23\}$. Then $X_v=\PSL_2(p)$ and $d=p{+}1$. For $p=17$ or $23$, the group $\PSL_2(p)$ has no transitive permutation representation of degree $p$, and thus it cannot occur as a subgroup of $\S_p$. Therefore, $p=7$ or $11$, and $G$, $X$ and $X_{\{u,v\}}$ are listed in Table \ref{PSL2q-graph}.
In fact, $X_{uv}$ and $X_{\{u,v\}}$ are the normalizers of some Syolw $p$-subgroup
in $\PSL_2(p)$ and $\S_p$, respectively.
(Note that, for $p=7$, the resulting graph is the point-plane non-incidence graph of $\PG(3,2)$.)
\begin{table}[ht]
\[
\begin{array}{|l|l|l|l|l|l|}
\hline
G&X & X_{\{u,v\}}& X_v& d& \mbox{Remark} \\ \hline
\S_p&\S_p & \ZZ_p{:}\ZZ_{p-1} & \PSL_2(p)& p{+}1& p\in \{7,11\}, \Ga\mbox{ bipartite}\\ \hline
\M_{11}&\M_{11} & 3^2{:}\Q_8.2& \M_{10}& 10 &\K_{11}\\ \hline
\J_1&\J_1& \ZZ_{11}{:}\ZZ_{10}& \PSL_2(11)& 12 &\\ \hline
\J_3.2&\J_3.2& \ZZ_{19}{:}\ZZ_{18}& \PSL_2(19)& 20 &\Ga\mbox{ bipartite}\\ \hline
\ON.2&\ON.2& \ZZ_{31}{:}\ZZ_{30}& \PSL_2(31)& 32 &\Ga\mbox{ bipartite}\\ \hline
\B&\B& \ZZ_{19}{:}\ZZ_{18}{\times} \ZZ_2& \PGL_2(19)& 20 & X_v<\Th<\B\\
 & & \ZZ_{23}{:}\ZZ_{11}{\times} \ZZ_2& \PSL_2(23)& 24 & X_v<\Fi_{23}<B\\ \hline
 \M&\M& \ZZ_{41}{:}\ZZ_{40} & \PSL_2(41)& 42 & \mbox{see \cite{Monster41} for } X_v \\ \hline
 \PSL_2(19)& \PSL_2(19)& \D_{20}& \PSL_2(5)& 6& \\ \hline
X, X.2 &  \PGL_2(9)& \D_{20}& \PSL_2(5)& 6& \K_{6,6}\\ \hline
X, X.2 &  \M_{10}& \ZZ_5{:}\ZZ_4& \PSL_2(5)& 6&  \K_{6,6}\\ \hline
 \PGL_2(11) &  \PGL_2(11)& \D_{20}& \PSL_2(5)& 6& \Ga\mbox{ bipartite} \\ \hline
 X, X.2  &  \PSL_3(r)& 3^2{:}\Q_8& \PSL_2(9)& 10&  r \mbox{ prime, \cite[Tables 8.3, 8.4]{Low}} \\

 &   &  && & r\equiv 4,16,31,34\,(\mod 45) \\ \hline

X, X.2  &  \PSU_3(r)& 3^2{:}\Q_8& \PSL_2(9)& 10&  r \mbox{ prime, \cite[Tables 8.5, 8.6]{Low}} \\

   &   &  && & r\equiv 11,14,29,41\,(\mod 45) \\ \hline

\end{array}
\]
{\caption{}
\label{PSL2q-graph}}
\end{table}

{\bf Case 2}.  Assume that $\soc(X)$ is a simple sporadic group. By \cite[Table 15]{s-arc}, with the restrictions (a) and (b), the only pairs $(G_0,H_0)$ are listed as follows{:}
$(\M_{11}, 3^2{:}\Q_8.2)$,
$(\J_1,\ZZ_{11}{:}\ZZ_{10})$, $(\J_1,\ZZ_{7}{:}\ZZ_{6})$, $(\J_3.2,\ZZ_{19}{:}\ZZ_{18})$, $(\J_4,\ZZ_{29}{:}\ZZ_{28})$,
$(\ON.2,\ZZ_{31}{:}\ZZ_{30})$, $(\B,\ZZ_{19}{:}\ZZ_{18}{\times} \ZZ_2)$, $(\B,\ZZ_{23}{:}\ZZ_{11}{\times} \ZZ_2)$, $(\M,\ZZ_{41}{:}\ZZ_{40})$ and $(\M,\ZZ_{47}{:}\ZZ_{23}{\times} \ZZ_2)$. In particular,  $\O_p(H_0)$ is  a Sylow $p$-subgroup of $G_0$. This yields that $X_v^{[1]}=1$, and so $\soc(X_v)=\PSL_2(p^f)$.

Suppose that $(X,X_{\{u,v\}})$ is one of $(\J_1,\ZZ_{7}{:}\ZZ_{6})$, $(\J_4,\ZZ_{29}{:}\ZZ_{28})$ and $(\M,\ZZ_{47}{:}\ZZ_{23}{\times} \ZZ_2)$.
Then  $X_v= \PSL_2(p)$ for $p=7$, $29$ and $47$, respectively; however,
by the Altas \cite{Atlas} and \cite[Tables 5.6 and 5.11]{Wilson-book}, $X$ has no subgroup $\PSL_2(p)$, a contradiction.
Thus $G$, $X$ and $X_{\{u,v\}}$ are listed in Table \ref{PSL2q-graph}.
(Note that the Monster $\M$ has a maximal subgroup $\PSL_2(41)$ by \cite{Monster41}.)

{\bf Case 3}.  Assume that $\soc(X)$ is a simple  group of Lie type over a finite field of order $r^t$, where $r$ is a prime. We first show $r\ne p$.

 Suppose that $r=p$. Then, by (a), either $\O_p(H_0)$ is  abelian or $r=p=2$. For $r=p>2$, noting that $|H_0|$ has a divisor $q-1$, there does not exist  $H_0$ in \cite[Tables 16-20]{s-arc} such that $\O_p(H_0)$ is  abelian.
Thus we have $r=p=2$. Recalling that $p^f>4$ and $|H_0/\O_p(H_0)|$ is divisible by $2^f-1$, it follows from Lemma \ref{2-3} that $H_0/\O_p(H_0)$ is not a $\{2,3\}$-group. Checking those $H_0$ given in \cite[Tables 16-20]{s-arc}, we conclude that $(G_0,H_0)$ is one of the following pairs:
 \begin{itemize}
 \item[] $(\PSL_2(2^t), \ZZ_2^t{:}\ZZ_{2^t-1})$,  $(\PSL_3(2^t), [2^{3t}]{:}[{(2^t-1)^2\over (3,2^t-1)}].2)$;
     $(\PSU_3(2^t), [2^{3t}]{:}\ZZ_{2^{2t}-1\over (3,2^t+1)})$;\\
  $(\PSp_4(2^t).\ZZ_{2^{l+1}}, [2^{4t}]{:}\ZZ_{2^{t}-1}^2.\ZZ_{2^{l+1}})$, where $2^l$ is the $2$-part of $t$;\\
      $(\Sz(2^t), [2^{2t}]{:}\ZZ_{2^{t}-1})$;
 $({}^3\D_4(2), [2^{11}]{:}(\ZZ_7\times\S_3))$; $({}^2\F_4(2)', [2^9]{:}5{:}4)$.
\end{itemize}
First, the pair $(\Sz(2^t), [2^{2t}]{:}\ZZ_{2^{t}-1})$ is excluded as $\Sz(2^t)$ has no subgroup with  a section $\PSL_2(2^f)$.
For the last two pairs, we have $f=5$ and $4$ respectively, which yields that
$2^f-1$ is not a divisor of $|H_0|$, a contradiction.
For the  three pairs after the first one, we have  $t<f$, thus $G_0$ has no maximal subgroup with a section $\PSL_2(2^f)$, a contradiction.
Suppose finally that $(X,X_{\{u,v\}})=(\PSL_2(2^t), \ZZ_2^t{:}\ZZ_{2^t-1})$. Then $3\le f<t\le 2f+1$.
Noting that $2^f-1$ is a divisor of $2^t-1$, it follows that $f$ is a divisor of $t$, and so $t=2f$. Then $\O_2(X_{\{u,v\}})=2^{2f}$, yielding $|\O_2(X_v^{[1]})|=2^f$.
Thus $\O_2(X_v)\ne 1$ and $X_v$ has a simple section
$\PSL_2(2^f)$. Checking the  subgroups of $\PSL_2(2^{2f})$, refer to \cite[II.8.27]{Huppert}, we conclude that $\PSL_2(2^{2f})$ has no    subgroup isomorphic to $X_v$, a contradiction.

Assume that   $r\ne p$ in the following.

{\it Subcase 3.1}. We first deal with those pairs $(G_0,H_0)$ such that $H_0$  is included in some infinite families in \cite[Table 16-20]{s-arc}.
Note that $r\ne p$, and we consider only those $H_0$ having subgroups of index $2$.
It follows that either $H_0/\Fit(H_0)$ is a $\{2,3\}$-group or
$(G_0,H_0)=(\E_8(q'), \ZZ_{q'^8{\pm} q'^7{\mp} q'^5-q'^4{\mp} q'^3{\pm} q'{+}1}.\ZZ_{30})$, where $q'=r^t$.
 Suppose that $(G_0,H_0)=(\E_8(q'), \ZZ_{q'^8{\pm} q'^7{\mp} q'^5-q'^4{\mp} q'^3{\pm} q'{+}1}.\ZZ_{30})$. Then $q'^8{\pm} q'^7{\mp} q'^5-q'^4{\mp} q'^3{\pm} q'{+}1$ is divisible by some primitive prime divisor $s$ of $q'^{15}-1$ or of $q'^{30}-1$. Noting that $s\ge 17$, we know that $H_0$ has normal cyclic Sylow $s$-subgroup.
It follows  from (a) that $17\le p=s=q'^8{\pm} q'^7{\mp} q'^5-q'^4{\mp} q'^3{\pm} q'{+}1$.
In particular, $\O_p(H_0)=\ZZ_p$ and $f=1$.
By (b), $|H_0|$ is divisible by $p-1$, and then $30$ is divisible by $p-1$.
This implies that $30=p-1=q'^8{\pm} q'^7{\mp} q'^5-q'^4{\mp} q'^3{\pm} q'$, which is impossible. Therefore, $H_0/\Fit(H_0)$ is a $\{2,3\}$-group.

By (a), $\Fit(H_0)$ a $\{2,p\}$-group. Then $|H_0|$ has no prime divisor other than $2$, $3$ and $p$. Since $p^f-1$ is a divisor of $|H_0|$, by Lemma \ref{2-3},
we have $f<3$.
Recall that $(X_u^{\Ga(u)})_v\cong (X_v^{\Ga(v)})_u= p^f{:}{q-1\over (2,q-1)}{.}[o]$,
 and
$X_{uv}\lesssim (X_v^{\Ga(v)})_u{\times} (X_u^{\Ga(u)})_v$, where $o$ is a divisor of $(2,q-1)f$. It follows that $X_{uv}/\O_p(X_{uv})$ has an abelian Hall $2'$-subgroup.
Note that $X_{uv}\O_p(X_{\{u,v\}})/\O_p(X_{\{u,v\}})\cong X_{uv}/(\O_p(X_{\{u,v\}}\cap X_{uv})=X_{uv}/\O_p(X_{uv})$, and $|X_{\{u,v\}}:X_{uv}\O_p(X_{\{u,v\}})|\le 2$.
It follows that $X_{\{u,v\}}/\O_p(X_{\{u,v\}})$ has an abelian Hall $2'$-subgroup.
Thus, as a possible candidate for $X_{\{u,v\}}$, the quotient of $H_0$ over $\O_p(H_0)$ has abelian Hall $2'$-subgroups. In particular, $H_0/\O_p(H_0)$ has no section $A_4$.

Considering the restrictions on $H_0$, $r$ and $f$, we conclude that
$(G_0,H_0)$ can only be one of the  following pairs:
\begin{itemize}
\item[] $(\PSL_2(r^t), \ZZ_{r^t\pm 1\over (2,r^t-1)}{:}\ZZ_2)$, $(\PSL_3(r^t), [{(r^t-1)^2\over (3,r^t-1)}].\S_3)$, $(\PSU_3(r^t), [{(r^t+1)^2\over (3,r^t+1)}].\S_3)$;\\
     $(\PSp_4(2^t).\ZZ_{2^{l+1}}, \ZZ_{2^t\pm 1}^2.[2^{l+4}])$, $(\PSp_4(2^t).\ZZ_{2^{l+1}}, \ZZ_{2^{2t}+1}.[2^{l+3}])$, where   $t\ge 3$;\\
$(\Sz(2^t), \ZZ_{2^t-1}{:}\ZZ_2)$,
$(\Sz(2^t), \ZZ_{2^t\pm \sqrt{2^{t+1}}+1}{:}\ZZ_4)$, where   $t\ge 3$;; \\
$(\Ree(3^t), \ZZ_{3^t\pm \sqrt{3^{t+1}}+1}{.}\ZZ_6)$,
$(\Ree(3^t), \ZZ_{3^t+1}{.}\ZZ_6)$,  where   $t\ge 3$;; \\
$(\G_2(3^t).\ZZ_{2^{l+1}}, \ZZ_{3^t\pm 1}^2{.}[3\cdot2^{l+3}])$,
$(\G_2(3^t).\ZZ_{2^{l+1}}, \ZZ_{3^{2t}\pm3^t+1}{.}[3\cdot2^{l+2}])$, where $t\ge 2$;\\
$({}^3\D_4(r^t), \ZZ_{r^{4t}-r^{2t}+1}{:}\ZZ_4)$;
$({}^2\F_4(2^t), \ZZ_{2^{2t}\pm \sqrt{2^{3t+1}}+2^{t}\pm \sqrt{2^{t+1}}+1}{.}\ZZ_{12})$, where   $t\ge 3$;
\\ $(\F_4(2^t).\ZZ_{2^{l+1}}, \ZZ_{2^{4t}-2^{2t}+1}{.}[3\cdot2^{l+3}])$, where $t\ge 2$;
\end{itemize}
where the power $2^l$ appeared means the $2$-part of $t$.
Recall that $|\Fit(H_0):\O_p(H_0)|\le 2$ and $|H_0:\O_p(H_0)|$ is divisible by $p^f-1$. This allows us determine the values of $p^f$ and $r^t$. As an example, we only deal with  the second pair. Suppose  that $(G_0,H_0)=(\PSL_3(r^t), [{(r^t-1)^2\over (3,r^t-1)}].\S_3)$. Considering the structures of $\Fit(H_0)$ and $\O_p(H_0)$, either $(3,r^t-1)=1$, $p=r^t-1$ and $f\in \{1,2\}$, or $f=1$ and $p=r^t-1=3$. The latter implies that $\PSL_2(q)$ is soluble, which is not the case. Assume that the former case holds.
Then $|\S_3|$ is divisible by $r^t-1-1$ or $(r^t-1)^2-1$. Then the only possibility is that $(p^f,r^t)=(7,8)$.
The other pairs can be fixed out in a similar way, the details is omitted here.
Eventually, we conclude that $(G_0,H_0, p,f)$ is one of
$(\PSL_2(19), \D_{20}, 5,1)$, $(\PSL_3(8), 7^2{:}\S_3,7,1)$ and
$(\Sz(8),\ZZ_5{:}\ZZ_4,5,1)$.
By the Atlas \cite{Atlas}, neither $\PSL_3(8)$ nor $\Sz(8)$ has subgroup
with a section $\PSL_2(p)$.
Thus, in this case, $G$, $X$ and $X_{\{u,v\}}$ are given as in Table \ref{PSL2q-graph}.

{\it Subcase 3.2}. For the pairs $(G_0,H_0)$ not appearing in {\it Subcase 3.1}, we check the finite number of $H_0$ one by one. We observed that
either $p=2$, or  $H_0/\O_p(H_0)$ is a $\{2,3\}$-group.
Recall that $r\ne p$.

Suppose that $p=2$. Recalling that $q=2^f>4$, we have $f\ge 3$.
In particular, since  $|H_0|$ is divisible by $2^f-1$,
$H_0$ is not a $\{2,3\}$-group
 by Lemma \ref{2-3}. Then the only possibility is that
$G_0={}^2\F_4(2)'$ and $H_0=[2^9]{:}5{:}4$. Thus $|\O_2(H_0)|=2^9$, it follows from (a) that $f=4$ or $9$, and then $G_0$ has a section $\PSL_2(2^4)$ or $\PSL_2(2^9)$, which is impossible by checking the (maximal) subgroups of ${}^2\F_4(2)'$. Thus $p>2$, and $H_0/\O_p(H_0)$ is a $\{2,3\}$-group; in particular, by (a),
 $\O_p(H_0)=\ZZ_p^{if}$ for some $i\in \{1,2\}$.

Suppose that $H_0$ has a section  $\A_4$.
Then $H_0$  has no normal Sylow $3$-subgroup. Further, $H_0$ has no quotient $\A_4$ as $H_0$ has  a subgroup of index $2$.
If $(3,(q-1)f)=1$ then, by (b), we conclude that $p=3$ and $\O_p(H_0)$ is the unique Sylow $3$-subgroup of $H_0$, a contradiction. Thus  $3$ is a divisor of $(q-1)f$.
Check those $H_0$ in \cite[Table 16-20]{s-arc} which have a section  $\A_4$ and
do not appear  in {\it Subcase 3.1}.
Recalling $r\ne p>2$ and $\O_p(H_0)=\ZZ_p^{if}$, it follows that either $\O_p(H_0)=\ZZ_3^2$ or $(G_0,H_0)=(\F_4(2).4, \ZZ_7^2{:}(3\times \SL_2(3)).4)$. Since $3$ is a divisor of $(q-1)f$, we get $G_0=\F_4(2).4$
and $q=p^f=7$ or $7^2$.
By (b), for $q=7$ or $7^2$, the order of $H_0$ should be a divisor of
$72$ or $192$ respectively, which is impossible.

The above argument  allows us ignore many cases without further inspection.
Inspecting carefully the remaining pairs,  the possible candidates for $(X,X_{\{u,v\}})$ are  as follows{:}
\begin{itemize}
\item[] $(\PGL_2(9), \D_{20})$, $(\M_{10}, \ZZ_5{:}\ZZ_4)$, $(\PGL_2(11), \D_{20})$;
\item[] $(\PSL_3(r), 3^2{:}\Q_8)$, where $r\equiv 4,7\,(\mod 9)$;
\item[] $(\PSp_4(4).4, \ZZ_{17}{:}\ZZ_{16})$, $(\PSp_4(4).4, 5^2{:}[2^5])$;
\item[] $(\PSU_3(r), 3^2{:}\Q_8)$, where $5<r\equiv 2,5\,(\mod 9)$;
\item[] $(\PSU_3(2^t), 3^2{:}\Q_8)$, where $t$ is a prime no less than $5$;
\item[] $({}^2\F_4(2), \ZZ_{13}{:}\ZZ_{12})$.
\end{itemize}

 For the first three pairs, $G$, $X$ and $X_{\{u,v\}}$ are easily determined and given as in Table \ref{PSL2q-graph}.
 The pair $(\PSp_4(4).4, \ZZ_{17}{:}\ZZ_{16})$ is excluded as $\PSp_4(4).4$ has no subgroup $\PSL_2(17)$, and the pair $({}^2\F_4(2), \ZZ_{13}{:}\ZZ_{12})$ is excluded as ${}^2\F_4(2)$  has no subgroup $\PSL_2(13)$. Suppose that $X\cong \PSU_3(2^t)$ and $X_{\{u,v\}}\cong 3^2{:}\Q_8$. Then we have $X_v\cong \PSL_2(9)$; however,
 by \cite[Tables 8.3, 8.4]{Low}, $\PSU_3(2^t)$ has no subgroup $\PSL_2(9)$, a contradiction.
 Suppose that   $(X,X_{\{u,v\}})=(\PSp_4(4).4, 5^2{:}[2^5])$.
 Then $X_v$ contains a Sylow $5$-subgroup $P$ of $X$ and has a section $\PSL_2(5)$ or
 $\PSL_2(25)$. By the information for $\PSp_4(4).4$ given in the Atlas \cite{Atlas},
 we conclude that $X_v\le M\cong (\A_5\times \A_5){:}2^2<\PSp_4(4).2<\PSp_4(4).4$.
 Note that $X_{uv}=5^2{:}[2^4]$, which should be the normalizer of $P$ in $X_v$.
 Using GAP \cite{GAP}, computation  shows that $|\N_L(P)|\le 200$ for any maximal subgroup $L$ of $M$ with $P\le L$. It follows that $X_v=M\cong (\A_5\times \A_5){:}2^2$,
 yielding $d=|X_v:X_{uv}|=36\ne q+1$, a contradiction.

Let $(X,X_{\{u,v\}})=(\PSL_3(r), 3^2{:}\Q_8)$. Then $X_{uv}\cong 3^2{:}4$.
It is easily shown that $p=3$ and $X_v\cong \PSL_2(9)$.
Since $r\equiv 4,7\,(\mod 9)$, we know that $\PSL_3(r)$ has a Sylow $3$-subgroup  $\ZZ_3^2$.
By \cite[Tables 8.3, 8.4]{Low}, $\PSL_3(r)$ has a subgroup $\PSL_2(9)$
 if and only if $r\equiv 1,4\,(\mod 15)$. Thus, in this case, we have
$r\equiv 11,14,29,41\,(\mod 45)$. For a subgroup $\PSL_2(9)$ of $\PSL_3(r)$, taking a Sylow $3$-subgroup $Q$ of $\PSL_2(9)$, the normalizers of $Q$ in $\PSL_2(9)$ and $\PSL_3(r)$ are (isomorphic to) $3^2{:}4$ and $3^2{:}\Q_8$, respectively.
Then these two normalizers of $Q$ can serve as the roles of $X_{uv}$ and $X_{\{u,v\}}$, respectively.
 Thus $X$ and $X_{\{u,v\}}$ are  given as in Table \ref{PSL2q-graph}.
Noting  that $G=XG_{\{u,v\}}$, we have $G_{\{u,v\}}/X_{\{u,v\}}\cong G/X\lesssim\Out(\PSL_3(r))\cong \S_3$, and so
$G=X.[m]$ and $G_{\{u,v\}}=X_{\{u,v\}}.[m]$, where $m$ is a divisor of $6$.
Thus $|G_{uv}{:}X_{uv}|=m$, since $|G_v{:}G_{uv}|=10=|X_v{:}X_{uv}|$, we have
$|G_v{:}X_v|=m$. By \cite[Table 8.4]{Low}, $\N_{\Aut\PSL_3(r)}(X_v)=X_v.2$.
Since $X_v\unlhd G_v$, it follows that $m\le 2$. Thus $G=X$ or $X.2$, and if $G=X.2$ then $G_v=X_v.2\cong \PGL_2(9)$ and $G_{\{u,v\}}\cong 3^2{:}\Q_8.2$.
The  pair $(\PSU_3(r), 3^2{:}\Q_8)$ is similarly dealt with.
\qed

\begin{lemma}\label{s3}
If  Lemma {\rm \ref{Xu-al-1}~(s3)} holds then $G$, $X$, $X_{\{u,v\}}$ and $X_v$ are listed in Table {\rm \ref{PSU3q-graph}}.
\end{lemma}
\begin{table}[ht]
\[
\begin{array}{|l|l|l|l|l|l|}
\hline
G&X & X_{\{u,v\}}& X_v& d& \mbox{Remark} \\ \hline
\HS.2&\HS.2 & [5^3]{:}[2^5] & \PSU_3(5){:}2& 126& \Ga \mbox{ bipartite}\\ \hline
\Ru&\Ru & [5^3]{:}[2^5] & \PSU_3(5){:}2& 126&  \\ \hline
\end{array}
\]
{\caption{}
\label{PSU3q-graph}}
\end{table}
\proof
Let $X_v^{\Ga(v)}= \PSU_3(q).[o]$ and $q=p^f>2$, where $p$ is a prime and $o\div 2(3,q{+}1)f$. Then
$(X_v^{\Ga(v)})_u=p^{f{+}2f}{:}{q^2-1\over (3,q{+}1)}.[o]$, and  $X_{uv}^{[1]}=1$ by Theorem \ref{double-star}. Thus
$|\O_p(X_{\{u,v\}})|=p^{3f}.a$, $p^{4f}.a$, $p^{5f}.a$ or $p^{6f}.a$, where $a$ is a divisor of $(2,p)$.
Moreover, $\O_p(X_{\{u,v\}})$ is nonabelian, and
$X_{\{u,v\}}/\O_p(X_{\{u,v\}})$ has a
subgroup $\ZZ_{(q^2-1)\over (3,q{+}1)}$. We next determine which pair $(G_0,H_0)$ in \cite[Tables 14-20]{s-arc} is a possible candidate for $(X,X_{\{u,v\}})$.
Note that
we may ignore those $H_0$ which either has no subgroup of index $2$ or has abelian maximal normal $p$-subgroup. In particular, $\soc(X)$ is not an alternating group.

{\bf Case 1}. Let  $(G_0,H_0)$ be a pair with $H_0$ included in some infinite families given in \cite[Table 16-20]{s-arc}. Since
$\O_p(X_{\{u,v\}})$ is nonabelian, we conclude that $(X,\O_p(X_{\{u,v\}}))$ is one of the following pairs:
\begin{itemize}
\item[] $(\PSL_3(p^t).2, [p^{3t}])$, $(\PGL_3(p^t).2,  [p^{3t}])$ (with $p=2$),\\ $(\PSU_3(p^t),  [p^{3t}])$, $(\PSp_4(p^t).\ZZ_{2^{l+1}},  [p^{4t}])$ (with $p=2$), \\ $(\Sz(p^t),  [p^{2t}])$, $(\Ree(p^t),  [p^{3t}])$ and $(\G_2(p^t).\ZZ_{2^{l+1}},  [p^{6t}])$,
\end{itemize}
 where $2^{l}$ is the $2$-part of $t$.
Check  the maximal subgroups of $\PSp_4(p^t).\ZZ_{2^{l+1}}$,  $\Sz(p^t)$ and  $\Ree(p^t)$, refer to \cite[Table 8.14]{Low}, \cite[Theorem 9]{Suzuki} and \cite[Theorem C]{Kleidman-2}, respectively.
We conclude that none of $\PSp_4(t^f).\ZZ_{2^{l+1}}$,  $\Sz(p^t)$ and  $\Ree(p^t)$ has maximal subgroups with a simple section $\PSU_3(q)$, and they are excluded.
For the first three and the last pairs, $|X/\O_p(X_{\{u,v\}})|$ is a divisor of $2(p^t-1)^2$, and $\O_p(X_{\{u,v\}})=[p^{3t}]$ or $[p^{6t}]$. Clearly, $t\le 2f$.

Suppose that $t=2f$. Then $\soc(X)=\PSL_3(q^2)$ or $\PSU_3(q^2)$, and $\O_p(X_{\{u,v\}})=[q^6]$. It follows that $\O_p(X_v^{[1]})=[q^3]$.
Thus $\O_p(X_v)\ne 1$ and $X_v$   has an almost simple quotient $\PSU_3(q).[o]$.
Checking Tables 8.3 and 8.5 given in \cite{Low}, we conclude that $X$ has no maximal subgroup
containing $X_v$, a contradiction.
If $t=f$ then we have $(X,\O_p(X_{\{u,v\}}))=(\G_2(p^t).\ZZ_{2^{l+1}},  [q^6])$, and
we get a similar contradiction by checking the maximal subgroups of $\G_2(p^t).\ZZ_{2^{l+1}}$.

Suppose that  $f\ne t<2f$. Then $f>1$.
Recalling that $X_{\{u,v\}}/\O_p(X_{\{u,v\}})$ has a
subgroup $\ZZ_{(q^2-1)\over (3,q{+}1)}$, we know that $p^{2f}-1$ is a divisor of $2(3,q+1)(p^t-1)^2$.
If $p^{2f}-1$ has a primitive prime divisor  say $s$, then $s\ge 2f+1\ge 5$, and $s$
is not a divisor of $2(3,q+1)(p^t-1)^2$, a contradiction.
It follows from Zsigmondy's theorem that $2f=6$ and $p=2$, and so $t=1$ or $2$.
Then $7$ is a divisor of $p^{2f}-1$ but not a divisor of $2(3,q+1)(p^t-1)^2$, a contradiction.

{\bf Case 2}.
Let  $(G_0,H_0)$ be one of the  pairs  in \cite[Table 15-20]{s-arc} which is not considered in Case 1.
Assume that $X_{\{u,v\}}/\O_p(X_{\{u,v\}})$ is a $\{2,3\}$-group.
Then $p^{2f}-1$ has no prime divisor other than $2$ and $3$.
It follows that $f=1$, and so $p=q>2$. Calculation shows that $p\in \{3,5,7\}$.
For $q=p=3$, it is easily shown that $X_{\{u,v\}}/\O_p(X_{\{u,v\}})$ is a $2$-group.
These observations yield that
 either $q=p=3$ and $X_{\{u,v\}}/\O_p(X_{\{u,v\}})$ is a $2$-group,
or  $X_{\{u,v\}}$ is not a $\{2,3\}$-group.

Recalling that $X_{\{u,v\}}/\O_p(X_{\{u,v\}})$
has a subgroup $\ZZ_{(q^2-1)\over (3,q{+}1)}$ and
$|\O_p(X_{\{u,v\}})|=p^{if}.a$ for $3\le i\le 6$,
it follows that $(X,X_{\{u,v\}})$ is one of the following pairs:
\begin{itemize}
\item[] $(\HS.2, [5^3]{:}[2^5])$, $(\Ru, [5^3]{:}[2^5])$,  $(\McL, [5^3]{:}3{:}8)$, $(\Co_2, [5^3]{:}4\S_4)$,\\ $(\Th, [5^3]{:}4\S_4)$, $(\J_4, [11^3]{:}(5\times 2\S_4))$.
\end{itemize}
Then $q=p\in \{5,11\}$ and $X_v^{[1]}=1$.
In particular, $\soc(X_v)=\PSU_3(p)$, and $X_{\{u,v\}}$ is the normalizer
$\N_X(P)$ of some Sylow $p$-subgroup $P$ of $X$.
Thus $X_{uv}=X_v\cap X_{\{u,v\}}\le \N_{X_v}(P)$.
For the pairs $(\HS.2, [5^3]{:}[2^5])$ and $(\Ru, [5^3]{:}[2^5])$,
by the Atlas \cite{Atlas},
$X_{\{u,v\}}$ is a normalizer of some Sylow $5$-subgroup,
which intersects a maximal subgroup $\PSU_3(5){:}2$ of $\soc(X)$ at $[5^3]{:}8{:}2$,
thus $G$, $X$ and $X_{\{u,v\}}$ are listed in Table \ref{PSU3q-graph}.
The other pairs are excluded as follows.

First, the group $\Th$ is excluded as it has no maximal subgroup with
a simple section $\PSU_3(5)$, refer to \cite[Table 5.8]{Wilson-book}.
For the pair $(\McL, [5^3]{:}3{:}8)$, by the Atlas \cite{Atlas},
we have $X_v=\PSU_3(5)$, and so $X_{uv}\le \N_{\PSU_3(5)}(P)=[5^3]{:}8$, which
contradicts that $|X_{\{u,v\}}:X_{uv}|=2$.
For the pair  $(\J_4, [11^3]{:}(5\times 2\S_4))$, by \cite[Table 5.8]{Wilson-book},
$X_v=\PSU_3(11).2$, yielding $X_{uv}\le \N_{X_v}(P)=[11^3]{:}(5\times 8{:}2)$,
we get a similar contradiction.
For the pair  $(X,X_{\{u,v\}})=(\Co_2, [5^3]{:}4\S_4)$, by the Atlas \cite{Atlas},  $X_v<\HS.2<\Co_2$. Checking the maximal subgroups
of $\HS.2$, we have $X_v=\PSU_3(5)$ or $X_v=\PSU_3(5){:}2$. It follows that
$X_{uv}\le \N_{X_v}(P)=[5^3]{:}8$ or $[5^3]{:}[2^5]$, and then $|X_{\{u,v\}}:X_{uv}|\ne 2$, a contradiction.
\qed

\vskip 20pt

\section{Graphs with soluble vertex-stabilizers}\label{sect-soluble}

Let $G$, $T$, $X$ and $\Ga=(V,E)$ be as in  Hypothesis \ref{hypo}.
We now deal with the case that $G_v$ is soluble. First, the following lemma says  that   $\Ga$ is not a complete bipartite graph if $G_v^{\Ga(v)}$ is soluble.

\begin{lemma}\label{complete-bi}
Assume that $\Ga\cong \K_{d,d}$. Then  $T\cong \A_6$,  $d=6$, $T_v=\PSL_2(5)$ and $T_{uv}\cong \D_{10}$. In particular, $X_{uv}$ is nonabelian.
\end{lemma}
\proof
Let $G^{+}$ be the subgroup of $G$ fixing the bipartition of $\Ga$.
Then $G_v\le G^{+}$, and $G_v$ is $2$-transitive on the partite set
which does not contain $v$. Thus $G^{+}$ acts $2$-transitively on
each partite set, and these two actions are not equivalent.
Check the almost simple $2$-transitive groups, refer to
\cite[Table 7.4]{Cameron}. We conclude that  $T\cong \A_6$ or $\M_{12}$,
and $T_v\cong \A_5$ or $\M_{11}$, respectively. Since $T_{uv}$ is soluble,
the lemma follows.
\qed

By Lemma \ref{complete-bi} and \cite{2-arc-edge-p}, we have the following result.

\begin{lemma}\label{abelian-arc-stab}
Assume that $X_{uv}$ is abelian. Then either $T\cong \PSL_2(q)$ and $\Ga\cong \K_{q{+}1}$, or $T=\Sz(2^{2m{+}1})$, $T_{\{u,v\}}\cong \D_{2(2^{2m{+}1}-1)}$, $T_v\cong \ZZ_2^{2m{+}1}{:}\ZZ_{2^{2m{+}1}-1}$, $\Aut\Ga=\Aut(\Sz(q))$ and $\Ga$ is $(T,2)$-arc-transitive.
\end{lemma}

\begin{lemma}\label{(a1)-(a2)}
Assume that {\rm (a1)} or {\rm (a2)} of Lemma {\rm \ref{stab-soluble-1}} holds,
and $X_{uv}$ is nonabelian. Then   one of the following holds{:}
\begin{enumerate}
\item[(1)] $G=X$ or $X.2$, $X=\M_{10}$, $X_{\{u,v\}}\cong \ZZ_8{:}\ZZ_2$,
 $X_v\cong 3^2{:}\Q_8$  and $\Ga\cong \K_{10}$.
\item[(2)] $G=X=\PSL_3(3).2$, $X_{\{u,v\}}\cong \GL_2(3){:}2$, $X_v\cong 3^2{:}\GL_2(3)$, and $\Ga$ is the point-line non-incidence graph of $\PG(2,3)$.
\end{enumerate}
\end{lemma}
\proof
{\bf Case 1}.  Assume that  Lemma {\rm \ref{stab-soluble-1}} (a1) holds.
Suppose that $(X_v^{\Ga(v)})_u=\Q_8$. Then $X_{uv}\lesssim\Q_8{\times}\Q_8$.
This implies that $|X_{\{u,v\}}|$ is a divisor of
$2^7$ and divisible by $2^4$. Checking the Tables 14-20 in \cite{s-arc},
we have $X\cong \PSL_2(9).2=\M_{10}$ and $X_{\{u,v\}}\cong \ZZ_8{:}\ZZ_2$;
in this case, $X_v\cong 3^2{:}\Q_8$, and $d=9$.
Then part (1) of this lemma follows.

Suppose that $(X_v^{\Ga(v)})_u\ne \Q_8$. If
 $p=3$ and $(G_v^{\Ga(v)})_u=\Q_8$, then  $(X_v^{\Ga(v)})_u$ is abelian,
 it follows that $X_{uv}$ is abelian, a contradiction. Thus we have
  $\SL_2(3)\unlhd (G_v^{\Ga(v)})_u\le \GL_2(p)$, and $p\in \{3,5,7,11,23\}$. Then $(G_v^{\Ga(v)})_u\le \N_{\GL_2(p)}(\SL_2(3))=\ZZ_{p-1}\circ \GL_2(3)$.
Since $(X_v^{\Ga(v)})_u$ is  nonabelian and normal in $(G_v^{\Ga(v)})_u$, we have $\Q_8\unlhd (X_v^{\Ga(v)})_u$, and hence
$\SL_2(3)\unlhd (X_v^{\Ga(v)})_u$. Moreover, $|X_{\{u,v\}}|$ is a divisor of
$2^73^2(p-1)^2$ and divisible by $2^4$.
Let $M$ be an arbitrary normal abelian subgroup of $X_{\{u,v\}}$. Then $M\cap X_{uv}$ has index at most $2$ in $M$,
and  $(M\cap X_{uv})X_v^{[1]}/X_v^{[1]}$ is isomorphic to a normal subgroup of
$(X_v^{\Ga(v)})_u$. Thus $(M\cap X_{uv})X_v^{[1]}/X_v^{[1]}\lesssim \ZZ_{p-1}$. Since $M\cap X_v^{[1]}\unlhd X_v^{[1]}$ and $X_v^{[1]}$ is isomorphic to a normal subgroup of
$(X_v^{\Ga(v)})_u$, we have $M\cap X_v^{[1]}\lesssim \ZZ_{p-1}$.
Noting that $(M\cap X_{uv})X_v^{[1]}/X_v^{[1]}\cong M\cap X_{uv}/(M\cap X_v^{[1]})$, it follows that
 $|M\cap X_{uv}|$ is a divisor of $(p-1)^2$. Thus $|M|$ is a divisor of $2(p-1)^2$.

The above observations allow us to consider only the pairs
$(G_0,H_0)$ in \cite[Tables 14-20]{s-arc} which satisfy the following conditions{:}
\begin{enumerate}
\item[(c1)] $|H_0|$ is a divisor of
$2^73^2(p-1)^2$ and divisible by $2^4$; $H_0$ has a factor (a quotient of some subnormal subgroup) $\Q_8$; and $H_0$ has no element of order $3^2$, $5^2$ or $11^2$.
 \item[(c2)] If $M$ is a normal abelian subgroup of $H_0$ then $|M|$ is a divisor of $2(p-1)^2$;
     if $p\in \{7,11,23\}$, the order of $\O_{p-1\over 2}(H_0)$ is a divisor of   ${(p-1)^2\over 4}$.
 \end{enumerate}
 Checking the those $H_0$ which satisfy conditions (c1) and (c2), we conclude that the possible pairs $(X,X_{\{u,v\}})$ are listed as follows{:}
 \begin{enumerate}
 \item[]
 $(\M_{11},3^2{:}\Q_8.2)$, $(\M_{11},2\S_4)$, $(\M_{12}, [2^5].\S_3)$,  $(\M_{12},3^2{:}2\S_4)$,  $(\J_{2},[2^6]{:}(3{\times}\S_3))$, \\ $(\J_{3},[2^6]{:}(3{\times}\S_3))$, $(\Co_3,[2^9].3^2.\S_3))$,  $(\He.2,[2^8]{:}3^2.\D_8))$, $(\McL.2, [2^6]{:}\S_3^2)$,  \\ $(\PSL_3(3), 3^2{:}2\S_4)$, $(\PSL_3(3).2, 2\S_4{:}2)$,
 $(\PSL_3(4).2, 2^{2{+}4}.3.2)$,
 \\ $(\PGL_3(4){.}2, [2^6]{.}3{.}\S_3)$, $(\PSL_4(3){.}2, 2{.}\S_4^2{.}2)$, $(\PSL_5(2){.}2, [2^8]{.}\S_3^2{.}2)$, \\
$(\PSp_4(4){.}4, [2^8]{:}3{.}12)$,  $(\PSp_4(4).4, 5^2{:}[2^5])$,
  $(\PSp_6(2), [2^7]{:}\S_3^2)$,  \\ $(\PSp_6(3), [2^8]{:}3^3.S_3)$,
  $(\PSU_3(3), 4.\S_4)$, $(\PSU_4(2), 2.\A_4^2.2)$,  $(\PSU_4(3), 2.\A_4^2.4)$,\\
  $(\PSU_4(3).2, [2^5].\S_4)$, $(\POmega^{+}_8(3){.}\A_4, 10^2{:}4\A_4)$, \\  $(\G_2(2)', 4.\S_4)$,  $(\G_2(3), \SL_2(3)\circ\SL_2(3){:}2)$, $({}^2\F_4(2)', 5^2{:}4\A_4)$.
  \end{enumerate}
 Note these groups $X$ are included in the Atlas \cite{Atlas}. Inspecting the subgroups of $X$, only the pair $(\PSL_3(3).2, 2\S_4{:}2)$ gives a desired $X_v\cong 3^2{:}\GL_2(3)$. Then part (2) of this lemma follows.

{\bf Case 2}. Let $2_{+}^{1{+}4}{:}\ZZ_5\le (G_v^{\Ga(v)})_u\le 2_{+}^{1{+}4}.(\ZZ_5{:}\ZZ_4)$. Then $2_{+}^{1{+}4}\unlhd (X_v^{\Ga(v)})_u$, and so, $|X_{\{u,v\}}|$ is a divisor of
$2^{15}5^2$ and divisible by $2^6$. Further, if $M$ is a normal abelian subgroup of
 $X_{\{u,v\}}$ then a similar argument as in {\bf Case 1} yields that $|M|$ is a divisor of $2^5$. It is easily shown that $\O_2(X_{uv})\ne1 $, and hence $\O_2(X_{\{u,v\}})\ne1$. Checking the pairs $(G_0,H_0)$  in \cite[Tables 14-20]{s-arc}, either $\O_2(H_0)=1$ or $|H_0|$ has an odd prime divisor other than $5$. Thus, in this case, no desired pair  $(X,X_{\{u,v\}})$ exists.
 \qed

  We assume next that Lemma {\rm \ref{stab-soluble-1}} (a3) occurs. Thus  $(G_v^{\Ga(v)})_u\not\le \GL_1(p^f)$ and $(G_v^{\Ga(v)})_u\le\GammaL_1(p^f)$. Then $f>1$ and $(G_v^{\Ga(v)})_u\lesssim\ZZ_{p^f-1}{:}\ZZ_f$. Recalling  $X_{uv}\lesssim (X_u^{\Ga(u)})_v{\times}  (X_v^{\Ga(v)})_u\le (G_u^{\Ga(u)})_v{\times}  (G_v^{\Ga(v)})_u$, we have the following simple fact.

\begin{lemma}\label{5-6}
If {\rm (a3)} of Lemma {\rm \ref{stab-soluble-1}} occurs then $X_{\{u,v\}}$ has no section $\ZZ_t^3$, $\ZZ_r^5$ or $\ZZ_2^6$, where $t$ is a primitive prime divisor of $p^f-1$ and $r$ is an arbitrary odd prime.
\end{lemma}

\begin{lemma}\label{(a3)-pf=2^6}
Assume that   $X_{uv}$ is nonabelian and   {\rm (a3)} of Lemma {\rm \ref{stab-soluble-1}} occurs.
Then  $p^f\ne 2^6$.
\end{lemma}
\proof
Suppose that $p^f=2^6$. Then $X$ has order divisible by $2^6$, $X_{uv}\lesssim \ZZ_{63}{:}\ZZ_6{\times} \ZZ_{63}{:}\ZZ_6$, and thus $X_{\{u,v\}}$ has a normal Hall $2'$-subgroup and
$|X_{\{u,v\}}|$ is indivisible by $2^4$. Checking Tables 14-20 given in \cite{s-arc},  $(X,X_{\{u,v\}})$ is one of the following pairs:
\begin{enumerate}
\item[] $(\S_7,\ZZ_7{:}\ZZ_6)$, $(\M_{12}.2,3_{+}^{1{+}2}{:}\D_8)$, $(\PSL_2(2^6),\D_{126})$,
$(\PSL_2(5^3),\D_{126})$, \\ $(\PSL_2(7937),\D_{7938})$, $(\PSL_3(8),7^2{:}\S_3)$,
$(\Sz(8), \D_{14})$, $(\G_2(3).2, [3^6]{:}\D_8)$.
\end{enumerate}
The pair $(\PSL_2(2^6),\D_{126})$ yields that $X_v\cong 2^6{:}\ZZ_{63}$, and thus
$X_{uv}$ is abelian, this is not the case. The other pairs are easily excluded as none of them gives a desired $X_v$.
\qed

\begin{lemma}\label{(a3)-abelian-2'}
Assume that   $X_{uv}$ is nonabelian and  {\rm (a3)} of Lemma {\rm \ref{stab-soluble-1}} occurs.
Suppose that  $X_{uv}$ has a normal abelian Hall $2'$-subgroup.
Then    $G=X$ or $X.2$, $X=\M_{10}$, $X_{\{u,v\}}\cong \ZZ_8{:}\ZZ_2$,
 $X_v\cong 3^2{:}\Q_8$  and $\Ga\cong \K_{10}$.
\end{lemma}
\proof
Note that $X_{\{u,v\}}=X_{uv}.2$.
The unique Hall $2'$-subgroup of $X_{uv}$ is also the Hall $2'$-subgroup of $X_{\{u,v\}}$.
Checking Tables 14-20 given in \cite{s-arc}, we know that $(X,X_{\{u,v\}})$ is one of the following pairs{:}
\begin{enumerate}
\item[(i)]$(\PGL_2(7),\D_{16})$, $(\PSL_3(2).2,\D_{16})$, $(\PGL_2(9),\D_{16})$, $(\M_{10},\ZZ_8{:}\ZZ_2)$; $(\A_5,\D_{10})$,\\ $(\A_6, 3^2{:}\ZZ_4)$,
$(\M_{11}, 3^2{:}\Q_8.2)$, $(\J_1,\D_6{\times} \D_{10})$, $(\PGL_2(7), \D_{12})$, $(\PGL_2(9),\D_{20})$, \\ $(\M_{10},\ZZ_5{:}\ZZ_4)$,
$(\PGL_2(11), \D_{20})$,
$(\PSL_2(t^a), \D_{2(t^a{\pm} 1)\over (2,t-1)})$; $(\PSp_4(4).4,\ZZ_{17}{:}\ZZ_{16})$;
 \item[(ii)]
 $(\PSL_2(t^a), \ZZ_t^a{:}\ZZ_{t^a-1\over 2})$, $t$ is a prime $a\le 4$ and $t^a-1$ a power of $2$;\\ $(\PSL_3(t),\ZZ_3^2{:}\Q_8)$, $t$ is a prime with $t\equiv 4,7\,(\mod 9)$; \\ $(\PSU_3(t),\ZZ_3^2{:}\Q_8)$, $t$ is a prime with $t\equiv 2,5\,(\mod 9)$; \\ $(\PSU_3(2^a),\ZZ_3^2{:}\Q_8)$ with prime $a>3$;\\
 $(\PSp_4(2^{a}).\ZZ_{2^{b{+}1}}, \D_{2(q{\pm} 1)}^2{:}2.\ZZ_{2^{b{+}1}})$,
 $(\PSp_4(2^{a}).\ZZ_{2^{b{+}1}}, \ZZ_{2^{2a}{+}1}.4.\ZZ_{2^{b{+}1}})$, where
 $2^b\parallel a$;\\
  $(\Sz(2^{2a{+}1}), \D_{2(2^{2a{+}1}-1)})$, $(\Sz(2^{2a{+}1}), \ZZ_{2^{2a{+}1}{\pm} 2^{a{+}1}{+}1}{:}\ZZ_4)$;\\
  $({}^3\D_4(t^a), \ZZ_{t^{4a}-t^{2a}{+}1}{:}\ZZ_4)$, $t$ is a prime.
\end{enumerate}

The pair $(\M_{10},\ZZ_8{:}\ZZ_2)$ yields that $X_v\cong 3^2{:}\Q_8$ and $d=9$.
The third pair in (i) implies that $X_v\cong \ZZ_3^2{:}\ZZ_8$; however, $X_{uv}$ is abelian, which is not the case.  For $(\PSL_2(t^a), \D_{2(t^a{\pm} 1)\over (2,t^a-1)})$, checking the subgroups of $\PSL_2(t^a)$,
we have $t^a=p^f$ and $X_v\cong \ZZ_p^f{:}\ZZ_{p^f-1\over (2,p-1)}$, and then $X_{uv}$ is abelian, a contradction.
The other pairs in (i) are also excluded as $|X|$ is indivisible by $p^f$. (Note that $f>1$.)

Now we deal with the pairs in (ii).
Note that, for an odd prime $r$, the edge-stabilizer  $X_{\{u,v\}}$ has a unique Sylow $r$-subgroup $\O_r(X_{\{u,v\}})$. Then $\O_r(X_{\{u,v\}})$ is a Sylow subgroup of $X$ by Lemma \ref{sylow}. This implies that the unique Hall $2'$-subgroup of $X_{\{u,v\}}$, say $K$, is a Hall subgroup of $X$. Since $X_{\{u,v\}}=X_{uv}.2$, we have $K\le X_{uv}$. Note that $|X_v{:}X_{uv}|=d=p^f$  and $X_v$ is contained in a maximal subgroup of $X$.
We now check the maximal subgroups of $X$ which contain $K$, refer to \cite[II.8.27]{Huppert},  \cite[Tables 8.3-8.6, 8.14, 8.15]{Low} and  \cite{Kleidman-1,Suzuki}. Then one of the following occurs{:}
\begin{enumerate}
\item[(iii)] $X=\Sz(2^{2a{+}1})$ and $X_v\cong \ZZ_2^{2a{+}1}{:}\ZZ_{2^{2a{+}1}-1}$;
\item[(iv)] $X=\PSp_4(2^{a}).\ZZ_{2^{b{+}1}}$ and $X_v\lesssim\Sp_2(2^{2a}){:}2.\ZZ_{2^b}$;
\item[(v)] $X=\PSp_4(2^{a}).\ZZ_{2^{b{+}1}}$ and $X_v\lesssim\Sp_2(2^{a})\wr \S_2.\ZZ_{2^b}$.
\end{enumerate}
Item (iii) yields that $X_{uv}$ is abelian, which is not the case.
Item (iv) gives $X_{uv}=X_v$, a contradiction.
Suppose that (v) occurs, we have $X_v\cong (\ZZ_2^{a}{:}\ZZ_{2^{a}-1})^2{:}2.\ZZ_{2^b}$.
Then $1\ne \O_2(X_v)\le \O_2(G_v)$, and hence
 $d=|\O_2(G_v)|$ by Lemma \ref{O_r(G-alpha)}. Since $X_v$ is transitive on $\Ga(v)$,
it follows that $p^f=d=2^{2a}$. Thus $|X_{uv}|=(2^{a}-1)^22^{b{+}1}$, and
so $|X_{\{u,v\}}{:}X_{uv}|=8>2$, a contradiction.
\qed

\begin{corollary}\label{(a3)-f=2}
Assume that   $X_{uv}$ is nonabelian and    {\rm (a3)} of Lemma {\rm \ref{stab-soluble-1}} occurs.
If  $f=2$ then $G=X$ or $X.2$, $X=\M_{10}$, $X_{\{u,v\}}\cong \ZZ_8{:}\ZZ_2$,
 $X_v\cong 3^2{:}\Q_8$  and $\Ga\cong \K_{10}$.
\end{corollary}
\proof
Let  $f=2$. Then $(X_v^{\Ga(v)})_u\lesssim\ZZ_{p^2-1}.\ZZ_{2}$. Note that $X_{\{u,v\}}=X_{uv}.2$ and $X_{uv}\lesssim \ZZ_{p^2-1}.\ZZ_{2}{\times} \ZZ_{p^2-1}.\ZZ_{2}$. Then
Lemma \ref{(a3)-abelian-2'} is applicable, and the result follows.
\qed

\vskip 5pt

Let ${\rm \pi}_0(p^f-1)$ be the set of  primitive
primes of $p^f-1$.
By  Zsigmondy's theorem, if ${\rm \pi}_0(p^f-1)=\emptyset$ and $f>1$ then
$p^f=2^6$, or $f=2$ and $p=2^t-1$, where $t$ is a prime.
Thus, in view of  Lemma \ref{(a3)-pf=2^6} and Corollary \ref{(a3)-f=2}, we assume next that ${\rm \pi}_0(p^f-1)\ne \emptyset$.

\begin{lemma}\label{(a3)-F}
Assume that  $\pi{:}={\rm \pi}_0(p^f-1)\ne \emptyset$, $X_{uv}$ is nonabelian and  {\rm (a3)} of Lemma {\rm \ref{stab-soluble-1}} occurs.
Then  $f\ge 3$, and
\begin{itemize}
\item[(1)] $\pi\ne \pi(|X_{\{u,v\}}|)\setminus\{2\}$, $\min(\pi)\ge \max\{5,f{+}1\}$;
\item[(2)] $p\not\equiv{\pm} 1\,(\mod r)$ and  $\O_r(X_{\{u,v\}})\ne 1$ for each $r\in \pi$;
\item[(3)] $X_{\{u,v\}}$ has a unique {\rm(}nontrivial\rm{)} Hall $\pi$-subgroup, which is either cyclic or  a direct product of two cyclic subgroups.
    \end{itemize}
\end{lemma}
\proof
By the assumption and Lemma \ref{(primitive-di)}, we have
$(X_v^{\Ga(v)})_u\cong \ZZ_{m'}.\ZZ_{f\over e'}$, and $\emptyset\ne \pi={\rm \pi}_0(p^f-1)\subseteq {\rm \pi}(m')$. For $r\in {\rm \pi}$, since $p^{r-1}\equiv 1\,(\mod r)$, we have $f\le r-1$, and so $r\ge f{+}1$. In particular, $r\ge 5$ and
$p\not\equiv{\pm} 1\,(\mod r)$. Recall  that $X_{\{u,v\}}=X_{uv}.2$ and $X_{uv}\lesssim \ZZ_{m'}.\ZZ_{f\over e'}{\times} \ZZ_{m'}.\ZZ_{f\over e'}$. It follows that
$\O_r(X_{\{u,v\}})\ne 1$, and $\O_r(X_{\{u,v\}})$ is the unique Sylow $r$-subgroup of $X_{\{u,v\}}$. Clearly, $\O_r(X_{\{u,v\}})$ is either cyclic or  a direct product of two cyclic subgroups.
Then $X_{\{u,v\}}$ has a unique Hall $\pi$-subgroup $F$, which is either cyclic or  a direct product of two cyclic subgroups. Clearly, $F\ne 1$ and, by Lemma \ref{(a3)-abelian-2'}, $X_{\{u,v\}}$ has no normal abelian Hall $2'$-subgroup.
Then $\pi\ne \pi(|X_{\{u,v\}}|)\setminus\{2\}$, the lemma follows.
\qed

\vskip 5pt

Recall that $X_{\{u,v\}}$ has no section  $\ZZ_2^6$ or $\ZZ_3^5$, see Lemma \ref{5-6}. Combining with Lemma \ref{(a3)-F}, we  next check the pairs $(G_0,H_0)$ listed in   \cite[Tables 14-20]{s-arc}.

\begin{lemma}\label{(a3)-Lie}
Assume that  ${\rm \pi}_0(p^f-1)\ne \emptyset$, $X_{uv}$ is nonabelian and   {\rm (a3)} of Lemma {\rm \ref{stab-soluble-1}} occurs.
Then  $T=\soc(X)$ is not a simple group of Lie type.
\end{lemma}
\proof
Suppose that $T$ is  a simple group of Lie type  over a finite field of order $q'=t^a$, where $t$ is a prime.
Since $T\unlhd G$, we know that $T$ is transitive on the edge set of $\Ga$. Then $T_v^{\Ga(v)}\ne 1$. Noting that $T_v^{\Ga(v)}\unlhd G_v^{\Ga(v)}$, we have $\soc(G_v^{\Ga(v)})\le T_v^{\Ga(v)}$. In particular, $T_v$ is transitive on $\Ga(v)$, and so $|T_v|=p^f|T_{uv}|$. In view of this, noting that $T_v=T\cap X_v=T\cap G_v$ and $T_{\{u,v\}}=T\cap X_{\{u,v\}}=T\cap G_{\{u,v\}}$, we sometimes work on the triple $(T,T_v,T_{\{u,v\}})$ instead of $(X,X_v,X_{\{u,v\}})$.

 By Lemmas \ref{(a3)-abelian-2'} and  \ref{(a3)-F}, $X_{\{u,v\}}$ is not a $\{2,3\}$-group and has no normal abelian Hall $2'$-subgroup.
Assume that  $t\in   \pi_0(p^f-1)$. By Lemmas \ref{5-6} and \ref{(a3)-F},  $t\ge 5$,   $X_{\{u,v\}}$   has no section $\ZZ_t^3$ and $\O_t(X_{\{u,v\}})\ne 1$ is abelian. Checking the pairs $(G_0,H_0)$ listed in   \cite[Tables 16-20]{s-arc}, we have  $X=\PSL_2(t^2)$ and $X_{\{u,v\}}\cong \ZZ_t^2{:}\ZZ_{t^2-1\over 2}$. For this case, checking the subgroups of $\PSL_2(t^2)$, no desired $X_v$ arises, a contradiction. Therefore,   $t\not\in   \pi_0(p^f-1)$.

 By Lemma \ref{(a3)-F}, $\O_r(X_{\{u,v\}})\ne 1$  for each $r\in \pi_0(p^f-1)$. Recall that
$X_{\{u,v\}}$ is not a $\{2,3\}$-group and has a subgroup of index $2$. Checking the pairs $(G_0,H_0)$ listed in   \cite[Tables 16-20]{s-arc}, we conclude that  $\O_t(X_{\{u,v\}})=1$. Further, we observe that a desired $X_{\{u,v\}}$ if exists has the form of
$N.K$, where $N$ is an abelian subgroup of $T$ and either $K$ is a $\{2,3\}$-group
or $(X,K)=(\E_8(q'),\ZZ_{30})$.
For the case where $K\not\cong \ZZ_{30}$, by Lemma \ref{(primitive-di)},  $\pi_0(p^f-1)\subseteq \pi(|N|)$, and thus, by Lemma \ref{5-6}, $N$ has no subgroup $\ZZ_r^3$ for $r\in \pi_0(p^f-1)$.  With these restrictions, only   one of the following (i)-(iv) occurs.

(i) Either $X=\PSL_3(q')$ and $X_{\{u,v\}}\cong {1\over (3,q'-1)}\ZZ_{q'-1}^2.\S_3$ with $q'\ne 2, \,4$, or $X=\PSU_3(q')$ and $X_{\{u,v\}}\cong {1\over (3,q'{+}1)}\ZZ_{q'{+}1}^2.\S_3$. Then $|X_v|={3\over (3,q'{\mp}1)}p^f(q'{\mp} 1)^2$. Checking Tables 8.3-8.6 given in \cite{Low}, we have $X=\PSL_3(q')$ and $X_v\lesssim[q'^3]{:}{1\over (3,q'-1)}\ZZ_{q'-1}^2$. It follows that $p=t=3$, and $|\O_3(X_v)|=3^{f{+}1}=3d$, which contradicts Lemma \ref{O_r(G-alpha)}.

(ii) $T=\soc(X)=\POmega_8^{+}(q')$ and $T_{\{u,v\}}\cong \D_{2(q'^2{+}1)\over (2,q'-1)}^2.[2^2]$. In this case, noting that $|T_{\{u,v\}}{:}T_{uv}|\le 2$, we have $|T_v|=2^4p^f{(q'^2{+}1)^2\over (2,q'-1)^2}$ or $2^3p^f{(q'^2{+}1)^2\over (2,q'-1)^2}$. Let $M$ be a maximal subgroup of $T$ with $T_v\le M$.
By \cite{Kleidman-1}, since $|M|$ is divisible by $(q'^2{+}1)^2$,
we  have $M\cong \PSL_2(q'^2)^2.2^2$. It is easily shown that $\PSL_2(q'^2)^2.2^2$
does not have subgroups of order $2^4p^f{(q'^2{+}1)^2\over (2,q'-1)^2}$ or $2^3p^f{(q'^2{+}1)^2\over (2,q'-1)^2}$, a contradiction.

(iii) $(X,X_{\{u,v\}})$ is one of
$({}^2\F_4(2)',5^2{:}4\A_4)$ and  $({}^2\F_4(2),13{:}12)$.
For the first pair, we have $\pi_0(p^f-1)=\{5\}$ and, since $p^f$ is a divisor of
$|{}^2\F_4(2)'|$, we conclude that
$p^f=2^4$ or $3^4$. The second pair implies that
 $\pi_0(p^f-1)=\{13\}$, and then $p^f=2^{12}$ or $3^3$. By the Atlas \cite{Atlas}, $X$ has no maximal subgroup containing $X_{uv}$ as a subgroup of index divisible by $p^f$, a contradiction.

(iv) $T_{\{u,v\}}$ has a normal abelian subgroup $N$   listed as follows{:}
\[
\begin{array}{l|l|c|l}  \hline
T & N & |T_{\{u,v\}}{:}N|& \mbox{Remark} \\ \hline
\Ree(3^a)& \ZZ_{3^a{\pm} 3^{a{+}1\over 2}{+}1}& 6 & \mbox{odd } a\ge 3\\
& \ZZ_2{\times} \ZZ_{3^a{+}1\over 2}& 6 & \&\, X=T\\  \hline

\G_2(3^a)& \ZZ_{3^a{\pm} 1}^2& 12 & \mbox{odd } a\ge 2\\
& \ZZ_{3^{2a}{\pm} 3^a{+}1}& 6 &  \\  \hline

& \ZZ_{2^a{+}1}^2& 48 & \mbox{odd } a\ge 3\\
{}^2\F_4(2^a)& \ZZ_{2^a{\pm} 2^{a{+}1\over 2}{+}1}^2& 96 & \&\, X=T\\
& \ZZ_{2^{2a}{\pm} 2^{3a{+}1\over 2}{+}2^a{\pm} 2^{a{+}1\over 2} {+}1}& 12 & \&\,2^a{\pm} 2^{a{+}1\over 2}{+}1> 5\\  \hline

& \ZZ_{2^{2a}{\pm} 2^a{+}1}^2& 72 & \\
\F_4(2^a) & \ZZ_{2^{2a}{+}1}^2& 96 & a\ge 2\\
& \ZZ_{2^{4a}-2^{2a} {+}1}& 12 &  \\  \hline

\E_8(q')& \ZZ_{q'^4-q'^2 {+}1}^2& 288 & X=T\\
& \ZZ_{q'^8{\pm} q'^7{\mp} q'^5-q'^4{\mp} q'^3{\pm} q'{+}1}& 30 &   \\  \hline
\end{array}
\]
Let $M$ be a maximal subgroup of $T$ with $T_v\le M$. Then $|M|$ is divisible by $p^f|N|$. Check the maximal subgroups of $T$ of order divisible by $|N|$, refer to
\cite{Kleidman-2, LSe, Malle}. Then we may deduce a contradiction.
First, by \cite[Theorem C]{Kleidman-2},   we conclude that $\Ree(3^a)$ has no maximal subgroups of order divisible by $p^f|N|$. Similarly, by  \cite{Malle}, the group ${}^2\F_4(2^a)$ is excluded. We next deal with the remaining cases.

(1) Let $T=\G_2(3^a)$.
Suppose that  $|N|=3^{2a}{\pm} 3^a{+}1$.
By \cite[Theorems A and B]{Kleidman-2}, since $|M|$ is divisible by $3^{2a}{\pm} 3^a{+}1$, we have   $M\cong \SL_3(3^a){:}2$ or $\SU_3(3^a){:}2$.
By \cite[Tables 8.3-8.6]{Low}, we conclude that $T_v\lesssim \ZZ_{3^{2a}{\pm} 3^a{+}1}{:}[6]$, which is impossible.

Similarly, for $|N|=(3^a{\pm} 1)^2$, we have
$T_v\lesssim (\SL_2(3^a)\circ\SL_2(3^a)).2$, $\SL_3(3^a){:}2$ or $\SU_3(3^a){:}2$.
Since $|T_v|$ is divisible by ${1\over 2}|T_{\{u,v\}}|p^f=6p^f(3^a{\pm} 1)^2$,
checking the maximal subgroups of $\SL_2(3^a)$, $\SL_3(3^a)$ and $\SU_3(3^a)$, we have $p=3$ and $T_v\lesssim [3^{ba}]{:}\ZZ_{3^a-1}^2.2$ for $b=2$ or $3$.
Since $T_{uv}$ has order divisible by $3$, it follows that $\O_3(T_{uv})\ne 1$, which contradicts Lemma \ref{O_r(G-alpha)}.

(3) Let $T=\F_4(2^a)$. By \cite{LSS,LSe}, noting that $|M|$ is divisible by $p^f|N|$, we conclude that $M\cong \Sp_8(2^a)$ or $\POmega_8^+(2^a).\S_3$ with $|N|=(2^{2a}{+}1)^2$, or $M\cong c.\PSL_3(2^a)^2.c.2$ or $c.\PSU_3(2^a)^2.c.2$ with $|N|=(2^{2a}{\pm} 2^a{+}1)^2$, where $c=(3,2^a\pm 1)$. Then a contradiction follows from checking the maximal subgroups of $\Sp_8(2^a)$, $\POmega_8^+(2^a)$, $\PSL_3(2^a)$ and $\PSU_3(2^a)$, refer to \cite[Tables 8.3-8.6, 8.48-8.50]{Low}.

(4) Let $T=\E_8(q')$. Then  $|N|=(q'^4-q'^2 {+}1)^2$
 and $M\cong \PSU_3(q'^2)^2.8$. For this case, checking the maximal subgroups $\PSU_3(q'^2)$, we get a contradiction.
\qed

\begin{lemma}\label{(a3)-3}
Assume that  ${\rm \pi}_0(p^f-1)\ne \emptyset$, $X_{uv}$ is nonabelian and   {\rm (a3)} of Lemma {\rm \ref{stab-soluble-1}} occurs.
Then    $G=X=\J_1$, $X_{\{u,v\}}\cong \ZZ_7{:}\ZZ_6$, $X_v\cong \ZZ_2^3{:}\ZZ_7{:}\ZZ_3$ and $d=8$.
\end{lemma}
\proof
By Lemma \ref{(a3)-Lie}, $T=\soc(X)$ is either an alternating  group or a sporadic simple group. Note that $X_{\{u,v\}}$ is not a $\{2,3\}$-group and has no normal abelian Hall $2'$-subgroup.

Assume that $T$ is an alternating  group. Then, by \cite[Table 14]{s-arc}, either $X=\A_r$ and $X_{\{u,v\}}\cong \ZZ_r{:}\ZZ_{r-1\over 2}$ for $r\not\in\{7,11,17,23\}$, or $X=\S_r$ and $X_{\{u,v\}}\cong \ZZ_r{:}\ZZ_{r-1}$ for $r\in\{7,11,17,23\}$. For these two cases, $X_v$ is a transitive subgroup of
$\S_r$ in the natural action of $\S_r$. Then either $X_v$ is almost simple or
 $X_v\lesssim \ZZ_r{:}\ZZ_{r-1}$ (refer to \cite[page 99, Corollary 3.5B]{Dixon}),
 a contradiction.

Assume that $T$ is a  sporadic simple group, and let $r\in \pi_0(p^f-1)$. Then  $(X,X_{\{u,v\}}, r)$ is one of the following triples{:}
\begin{enumerate}
\item[]
$(\J_1, \ZZ_7{:}\ZZ_6,7)$, $(\J_1, \ZZ_{11}{:}\ZZ_{10},11)$, $(\J_1, \ZZ_{19}{:}\ZZ_{6},19)$, $(\J_2, \ZZ_{5}^2{:}\D_{12},5)$, \\ $(\J_3.2, \ZZ_{19}{:}\ZZ_{18},19)$, $(\J_4, \ZZ_{29}{:}\ZZ_{28},29)$, $(\J_4, \ZZ_{37}{:}\ZZ_{12},37)$, $(\J_4, \ZZ_{43}{:}\ZZ_{14},43)$, \\
$(\ON.2, \ZZ_{31}{:}\ZZ_{30},31)$,    $(\He, \ZZ_5^2{:}4\A_4, 5)$,   $(\Co_1, \ZZ_{7}^2{:}(3{\times} 2\A_4),7)$,\\ $(\Ly, \ZZ_{37}{:}\ZZ_{18},37)$,  $(\Ly, \ZZ_{67}{:}\ZZ_{22},67)$, $(\Fi_{24}', \ZZ_{29}{:}\ZZ_{14},29)$,  \\  $(\B, \ZZ_{13}{:}\ZZ_{12}{\times} \S_4,13)$,  $(\B, \ZZ_{19}{:}\ZZ_{18}{\times}\ZZ_2,19)$,  $(\B, \ZZ_{23}{:}\ZZ_{11}{\times} 2,23)$,  \\
 $(\M, \ZZ_{23}{:}\ZZ_{11}{\times} \S_4,23)$,  $(\M, (\ZZ_{29}{:}\ZZ_{14}{\times} 3).2,29)$, $(\M, \ZZ_{31}{:}\ZZ_{15}{\times} \S_3,31)$, \\ $(\M, \ZZ_{41}{:}\ZZ_{40},41)$, $(\M, \ZZ_{47}{:}\ZZ_{23}{\times} 2,47)$.
\end{enumerate}
Recall that $p^f$ is a divisor of $|X|$ and $r$  is a primitive prime divisor of $p^f-1$. Searching all possible pairs $(p^f,r)$, we get the following table{:}
\[
\begin{array}{l|l|l|l|l|l|l|l}
X&\J_1 &\J_2 &\J_4&\Co_1&\ON.2&\He&\B\\ \hline

|X_{\{u,v\}}|& 2{\cdot} 3{\cdot} 7&   2^2{\cdot} 3{\cdot} 5^2 & 2{\cdot} 7{\cdot} 43&  2^3{\cdot} 3^2{\cdot} 7^2& 2{\cdot} 3{\cdot} 5{\cdot} 31& 2^4{\cdot}3{\cdot}5^2& 2^5{\cdot} 3^4{\cdot} 13\\ \hline

r& 7&   5&  43& 7&  31& 5& 13\\ \hline
p^f & 2^3&  2^4& 2^{14}&  2^3,\, 3^6& 2^5&2^4& 3^3,5^4,2^{12}\\  \hline

{p^f-1}\div |G_{uv}|& \checkmark& \checkmark& {\times} & \checkmark,\, {\times}&\checkmark& \checkmark& \checkmark,\checkmark,{\times}  \\  \hline   \hline

X&\B& \B&  \M& \M& \M&\M&\M  \\  \hline

|X_{\{u,v\}}|&2^2{\cdot} 3^2{\cdot} 19& 2{\cdot} 11{\cdot} 23&  2^3{\cdot} 3{\cdot} 11{\cdot} 23&
 2^2{\cdot} 3{\cdot} 7{\cdot} 29&  2{\cdot} 3^2{\cdot} 5{\cdot} 31& 2^3{\cdot} 5{\cdot} 41&2{\cdot} 23{\cdot} 47 \\  \hline

r&19& 23& 23&  29& 31&41&47 \\  \hline

p^f&2^{18} &  2^{11}, 3^{11}& 2^{11}, 3^{11}&2^{28}&2^5,5^3&2^{20}, 3^8&2^{23} \\ \hline

{p^f-1}\div |G_{uv}|&{\times}& {\times},\,{\times}&{\times},\,{\times}&{\times}&\checkmark,{\times}&{\times},\,{\times}&{\times}  \\

\end{array}
\]

 Recalling that $G_{\{u,v\}}=X_{\{u,v\}}.(G/X)$, we have
 $2|G_{uv}|=|G_{\{u,v\}}|=|X_{\{u,v\}}||G{:}X|=2|X_{uv}||G{:}X|$, and so
 $|G_{uv}|=|X_{uv}||G{:}X|$. Since $G_v$ is $2$-transitive on $\Ga(v)$, we know that $(p^f-1)$ is a divisor of $|G_{uv}|=|X_{uv}||G{:}X|$. It follows that $(X,X_{\{u,v\}}, r, p^f)$ is one of
 \begin{itemize}
\item[]  $(\J_1, \ZZ_7{:}\ZZ_6,7, 2^3)$, $(\J_2, \ZZ_{5}^2{:}\D_{12},5, 2^4)$, $(\Co_1, \ZZ_{7}^2{:}(3{\times} 2\A_4),7, 2^3)$,
 \item[] $(\ON.2, \ZZ_{31}{:}\ZZ_{30},31, 2^5)$, $(\He, \ZZ_5^2{:}4\A_4, 5,2^4)$,  $(\B, \ZZ_{13}{:}\ZZ_{12}{\times} \S_4,13, 3^3)$,
 \item[] $(\B, \ZZ_{13}{:}\ZZ_{12}{\times} \S_4,13, 5^4)$, $(\M, \ZZ_{31}{:}\ZZ_{15}{\times} \S_3,31, 2^5)$.
\end{itemize}

For $(\Co_1, \ZZ_{7}^2{:}(3{\times} 2\A_4),7, 2^3)$, we have $X_{uv}\lesssim \GammaL_1(2^3){\times} \GammaL_1(2^3)$, yielding that $X_{uv}$ has odd order, a contradiction. Similarly, for $(\B, \ZZ_{13}{:}\ZZ_{12}{\times} \S_4,13, 3^3)$, the order of  $X_{uv}$ is indivisible by $2^4$,   a contradiction; for $(\M, \ZZ_{31}{:}\ZZ_{15}{\times} \S_3,31, 2^5)$, the order of  $X_{uv}$ is indivisible by $3$, also  a contradiction.
For $(\He, \ZZ_5^2{:}4\A_4, 5,2^4)$, the order of  $X_{uv}$ is divisible by
$2^3\cdot 3\cdot 5^2$ and, since $p^f=2^4$, the order of $X_u$ is  divisible by
$2^7\cdot 3\cdot 5^2$; however, $\He$ has no soluble subgroups   of order divisible by $2^7\cdot 3\cdot 5^2$, a contradiction.
Similarly, the quadruple  $(\ON.2, \ZZ_{31}{:}\ZZ_{30},31, 2^5)$ is excluded as $\ON.2$ has no soluble subgroups   with order divisible by $2^5\cdot 31$. (Note that $G_v$ is soluble.) By the Altas \cite{Atlas}, $\J_2$ has no subgroups with order divisible by $2^4\cdot 5^2$, and then $(\J_2, \ZZ_{5}^2{:}\D_{12},5, 2^4)$ is excluded.
By the Altas \cite{Atlas} and \cite[Theorem 2.1]{Wilson-99}, $\B$ has no  subgroups
 with order divisible by $3^2\cdot 5^4\cdot 13$, and then $(\B, \ZZ_{13}{:}\ZZ_{12}{\times} \S_4,13, 5^4)$ is excluded. Then $(\J_1, \ZZ_7{:}\ZZ_6,7, 2^3)$ is left, which gives   $X_v\cong \ZZ_2^3{:}\ZZ_7{:}\ZZ_3$, $d=p^f=8$ and $G=X$.
\qed

\vskip 10pt

Finally, we summarize the argument for proving Theorem \ref{main-result} as follows.

\noindent{\it Proof of Theorem {\rm \ref{main-result}}}.
Clearly, each $(G,G_v, G_{\{u,v\}})$ in Table \ref{all-graphs}
gives a $G$-edge-primitive graph $\Cos(G,G_v, G_{\{u,v\}})$.
It is not difficult to check the $2$-arc-transitivity of $G$ on $\Cos(G,G_v, G_{\{u,v\}})$, we omit the details.

Now let $G$ and $\Ga=(V,E)$ satisfy the assumptions in Theorem \ref{main-result}.
If  $G_v^{\Ga(v)}$ is an almost simple $2$-transitive group then, by Lemmas \ref{Xu-al-1}, \ref{s4-s5}-\ref{PSL_2(q)}, the triple $(G,G_v, G_{\{u,v\}})$ is listed in Table \ref{all-graphs}. Assume that  $G_v^{\Ga(v)}$ is a soluble $2$-transitive group
of degree $d=p^f$, where $p$ is a prime. Then either  $G_v^{\Ga(v)}\le \GammaL_1(p^f)$, or $G_v^{\Ga(v)}$ has a normal subgroup $\SL_2(3)$ or $2_+^{1+4}$.
For the latter case, the triple $(G,G_v, G_{\{u,v\}})$ is known by Lemmas \ref{abelian-arc-stab} and \ref{(a1)-(a2)}.  Let $G_v^{\Ga(v)}\le \GammaL_1(p^f)$
and consider the primitive prime divisors of $p^f-1$.
If $p^f-1$ has no primitive prime divisor then, by  Lemmas \ref{abelian-arc-stab}, \ref{(a3)-pf=2^6} and Corollary \ref{(a3)-f=2}, $(G,G_v, G_{\{u,v\}})$ is listed in Table \ref{all-graphs}. If $p^f-1$ has   primitive prime divisors,  then $(G,G_v, G_{\{u,v\}})$ is known by Lemmas \ref{abelian-arc-stab} and  \ref{(a3)-3}.
This completes the proof.
\qed

\end{document}